\documentclass{amsart}

\usepackage{amsmath} 
\usepackage{amssymb}
\usepackage{mathrsfs}

\newtheorem{theorem}{Theorem}[section] 
\newtheorem{claim}[theorem]{Claim}

\theoremstyle{definition}
\newtheorem{definition}[theorem]{Definition}

\newtheorem{problem}[theorem]{Problem}
\newtheorem{observation}[theorem]{Observation} 

\newtheorem{explanation}[theorem]{Explanation}

\theoremstyle{remark}
\newtheorem{remark}[theorem]{Remark}
\newtheorem{notation}[theorem]{Notation}
\newtheorem{conclusion}[theorem]{Conclusion}

\newcommand{\rest}{{\restriction}}

\newcommand{\wilog}{{\rm without loss of generality}}

\newcommand{\then}{{\underline{then}}}
\newcommand{\when}{{\underline{when}}}

\newcommand{\Iff}{{\underline{iff}}}
\newcommand{\mn}{{\medskip\noindent}}
\newcommand{\sn}{{\smallskip\noindent}}

\newcommand{\cA}{{\mathcal A}}

\newcommand{\cB}{{\mathcal B}}

\newcommand{\gb}{{\mathfrak b}}

\newcommand{\ga}{{\mathfrak a}}
\newcommand{\cC}{{\mathcal C}}

\newcommand{\gS}{{\mathfrak S}}

\newcommand{\cP}{{\mathcal P}}

\newcommand{\gs}{{\mathfrak s}}

\newcommand{\cT}{{\mathcal T}}
 
\newcommand{\cU}{{\mathcal U}}
\newcommand{\cQ}{{\mathcal Q}}
\newcommand{\cW}{{\mathcal W}}

\newcount\skewfactor
\def\mathunderaccent#1#2 {\let\theaccent#1\skewfactor#2
\mathpalette\putaccentunder}
\def\putaccentunder#1#2{\oalign{$#1#2$\crcr\hidewidth
\vbox to.2ex{\hbox{$#1\skew\skewfactor\theaccent{}$}\vss}\hidewidth}}

\newenvironment{PROOF}[2][\proofname.]
   {\begin{proof}[#1]\renewcommand{\qedsymbol}{\textsquare$_{\rm #2}$}}
   {\end{proof}}

\begin{document}

\title {MAD saturated families and sane player}
\author {Saharon Shelah}
\address{Einstein Institute of Mathematics\\
Edmond J. Safra Campus, Givat Ram\\
The Hebrew University of Jerusalem\\
Jerusalem, 91904, Israel\\
 and \\
 Department of Mathematics\\
 Hill Center - Busch Campus \\ 
 Rutgers, The State University of New Jersey \\
 110 Frelinghuysen Road \\
 Piscataway, NJ 08854-8019 USA}
\email{shelah@math.huji.ac.il}
\urladdr{http://shelah.logic.at}
\thanks{Research supported by the United States-Israel Binational
Science Foundation (Grant No. 2006108).  Publication 935.
The author thanks Alice Leonhardt for the beautiful typing}

\date{June 18, 2010}

\begin{abstract}
We throw some light on the question: is there a MAD family
(= a maximal family of infinite subsets of $\Bbb N$, the intersection of any
two is finite) which is saturated (= completely separable i.e. any 
$X \subseteq \Bbb N$ is
included in a finite union of members of the family \underline{or} includes a
member (and even continuum many members) of the family).  
We prove that it is hard to prove the consistency of the negation:
\mn
\begin{enumerate}
\item[$(a)$]   if $2^{\aleph_0} < \aleph_\omega$, then there is such a family
\sn
\item[$(b)$]  if there is no such family then some situation
related to pcf holds whose consistency is large; and if ${\frak a}_* >
\aleph_1$ even unknown
\sn
\item[$(c)$]   if, e.g. there is no inner model with measurables \then \,
there is such a family.
\end{enumerate}
\end{abstract}

\maketitle
\numberwithin{equation}{section}
\setcounter{section}{-1}

\section{Introduction}

We try to throw some light on

\begin{problem}
\label{0z.2}  
Is there, provably in ZFC, a completely
separable MAD family ${\cA} \subseteq [\omega]^{\aleph_0}$, see
Definition \ref{0z.4}(1),(4).

Erd\"os-Shelah \cite{ErSh:19} investigates the ZFC-existence of
families ${\cA} \subseteq {\cP}(\omega)$ with separability
properties, continuing Hechler \cite{He71} which mostly uses MA; now
\ref{0z.2} is Problem A of \cite{ErSh:19}, pg.209, see
earlier Miller \cite{Mil37}, and see later
Goldstern-Judah-Shelah \cite{GJSh:399} on existence for larger
cardinals.  It seemed natural to prove the
consistency of a negative answer by CS iteration making the continuum
$\aleph_2$ but this had not worked out; the results here show this 
is impossible.

The celebrated matrix-tree theorem of Balcar-Pelant-Simon \cite{BPS},
Balcar-Simon \cite{BaSi89} is related to our starting point.
In Gruenhut-Shelah \cite{GhSh:E64} we try to generalize
it, hoping eventually to get applications, e.g. ``there is a subgroup
of ${}^\omega \Bbb Z$ which is reflexive (i.e. canonically isomorphic
to the dual of its dual)" and ``less", see Problem D7 of
\cite{EM02}, no success so far.
We then had tried to use such constructions to answer
\ref{0z.2} positively, but this does not work.
 Simon \cite{Si96} have proved (in ZFC), that
there is an infinite almost disjoint ${\cA} \subseteq
[\omega]^{\aleph_0}$ such that $B \subseteq \omega$ and
$(\exists^\infty A \in {\cA})[B \cap A \text{ infinite}] \Rightarrow
(\exists A \in {\cA})(A \subseteq B)$.  Shelah-Steprans
\cite{ShSr:931} try to continue it with dealing with Hilbert spaces.

Here ${\gs}$ and ideals (formally $J \in$ OB) are central.
Originally we have a unified proof using games between the MAD and the
SANE players but with some parameters for the
properties. As on the one hand it was claimed this is unreadable and
on the other hand we have a direct proof, which was presented (for
${\gs} < \ga_*$), in the Hebrew University and Rutgers, we
use the later one.  A minor price is that the proof in \S2 are saying -
repeat the earlier one with the following changes.  The major price is
that some information is lost: using smaller more complicated cardinal
invariants as well as some points in the proof which we hope will
serve other proofs (including covering all cases) so we shall
return to the main problem and relatives in 
\cite{Sh:F1047} which continue this work.

A related problem of Balcar and Simon is:  given a MAD family
${\cB}$ we look for such ${\cA}$ refining it, i.e. $(\forall B \in
\text{ id}^+_{\cA})(\exists A \in {\cA})(A \subseteq^* B)$.  At present there
is no difference between the two problems (i.e. in \ref{4d.3},
\ref{g.3}, \ref{g.21} we cover this too)
\end{problem}

Anyhow
\begin{conclusion}
\label{0z.3} 
1) If $2^{\aleph_0} < \aleph_\omega$
\then \, there is a saturated MAD family.

\noindent
2) Moreover in (1) for any dense $J_* \subseteq
   [\omega]^{\aleph_0}$ we can find such a family $\subseteq J_*$.
\end{conclusion}

We thank Shimoni Garti and the referee for helpful corrections.
\begin{definition}
\label{0z.4} 
1) We say ${\cA}$ is an AD (family) for $B$ 
when ${\cA} \subseteq [B]^{\aleph_0}$ is
 infinite, almost disjoint (i.e. $A_1 \ne A_2 \in {\cA}
 \Rightarrow A_1 \cap A_2$ finite).
We say ${\cA}$ is MAD for $B$ \when \, ${\cA}$ is AD for
$B$ and is $\subseteq$-maximal among such ${\cA}$'s.

\noindent
2) If $B = \omega$ we may omit it.

\noindent
3) For ${\cA} \subseteq [\omega]^{\aleph_0}$, id$_{\cA}$ is the
   ideal generated by ${\cA} \cup [\omega]^{< \aleph_0}$.

\noindent
4) A MAD family ${\cA}$ is saturated when: if $B \in
\text{ id}^+_{\cA}$ (see \ref{0z.6}(3)) then $B$ 
almost contains some  member of ${\cA}$ 
(equivalently: if $B \in \text{ id}^+_{\cA}$ 
\then \, $B$ almost contains continuum many members of $A$ because if
$B \in \text{ id}^+_{\cA}$ then there is an AD family $\cB \subseteq
[B]^{\aleph_0} \cap \text{ id}^+_{\cA}$ of cardinality $2^{\aleph_0}$).
\end{definition}

\begin{definition} 
\label{0z.3d}
\item[(1)]
Let ${\frak a}$ be the minimal cardinality of a MAD family
\item[(2)]
Let $ {\frak a}_* $ be the minimal $ \kappa $
such that there is a sequence 
$ \langle A_ \alpha : \alpha < \kappa + \omega \rangle $ 
of pairwise almost disjoint  (=with finite intersection)
infinite subsets of $ \omega $ satisfying:

there is no infinite set $B \subseteq \omega$ almost disjoint
to $A_ \alpha$ for $\alpha < \kappa$ but $B \cap A_{\kappa + n}$ 
is infinite for infinitely many $n$-s.
\end{definition}

\begin{observation} 
\label{0z.3f}
We have $ {\frak b} \le {\frak a}_* \le {\frak a}$. 
\end{observation}

\begin{remark}
\label{0z.5}  
1) Note that if there is a MAD family ${\cA} \subseteq
[\omega]^{\aleph_0}$ such that $B \in \text{ id}^+_{\cA}
\Rightarrow (\exists^{2^{\aleph_0}} A \in {\cA})$ ($B \cap A$ is
infinite), \underline{then} there is a MAD family ${\cA} \subseteq
[\omega]^{\aleph_0}$ such that $B \in \text{ id}^+_{\cA}
\Rightarrow (\exists^{2^{\aleph_0}} A \in {\cA})(A \subseteq B)$
equivalently $B \in \text{ id}^+_A \Rightarrow (\exists A \in 
{\cA})(A \subseteq B)$; just list our tasks and fulfil them by 
dividing each member of ${\cA}$ to two infinite sets to fulfil on task.

\noindent
2) So the four variants of ``there is ${\cA} \ldots$" in
\ref{0z.4}(4), \ref{0z.5}(1) are equivalent.
\end{remark}

\begin{notation}
\label{0z.6} 
1) For $A \subseteq \omega$ let
$A^{[\ell]}$ be $A$ if $\ell=1$ and $\omega \backslash A$ if $\ell=0$.

\noindent
2) For $J \subseteq [\omega]^{\aleph_0}$ let $J^\perp = \{B:B \in
   [\omega]^{\aleph_0}$ and $[A \in J \Rightarrow A \cap B \text{
   finite}]\}$ and also for $\bar A = \langle A_s:s \in S\rangle$ let
$\bar A^\perp = \{A_s:s \in S\}^\perp$.

\noindent
3) id$_{\cA}(B)$ is the ideal of ${\cP}(B)$ generated by 
$({\cA} \rest B) \cup [B]^{< \aleph_0}$ and id$^+_{\cA}(B) =
[B]^{\aleph_0} \backslash \text{ id}_{\cA}(B)$, for $\cA \rest B$ see
7) below; if $B = \omega$ we may omit it.

\noindent
4) Let $A \subseteq^* B$ means that $A \backslash B$ is finite.

\noindent
5) If ${\cC} \subseteq {\cP}(\omega)$ and $\eta \in 
{}^{\cC}2$ \then \, $I_{{\cC},\eta}(B)$ is $\{C \subseteq B:C
\subseteq^* A^{[\eta(A)]}$ for every $A \in {\cC}\}$; 
if $B=\omega$ we may omit it.

\noindent
6) In part 5), if $\nu$ is a function extending $\eta$ then let
$I_{{\cC},\nu} = I_{{\cC},\eta}$.

\noindent
7) For ${\cA} \subseteq {\cP}(B_2)$ and $B_1 \subseteq B_2$ let
${\cA} \rest B_1 = \{A \cap B_1:A \in {\cA}$ satisfies $A \cap
 B_1$ is infinite$\}$.
\end{notation}

\begin{definition}
\label{0z.8} 
1) Let OB $= \{I \subseteq [\omega]^{\aleph_0}:
I \cup [\omega]^{< \aleph_0}$ is an ideal of ${\cP}(\omega)\}$.

\noindent
2) For $A \subseteq \omega$ let ob$(A) = \{B:B \in
[\omega]^{\aleph_0}$ and $B \subseteq^* A\}$ so ob$(\omega) =
[\omega]^{\aleph_0}$.

\noindent
3) $\eta \perp \nu$ means $\neg(\eta \trianglelefteq \nu) \wedge
\neg(\nu \trianglelefteq \eta)$.

\noindent
4) We say ${\cA}$ is AD in $J \subseteq [\omega]^{\aleph_0}$
\when \, ${\cA}$ is AD and ${\cA} \subseteq J$.

\noindent
5)  We say ${\cA}$ is MAD in $J \subseteq [\omega]^{\aleph_0}$
\when \, ${\cA}$ is AD in $J$ and is $\subseteq$-maximal among
such $\cA$'s.

\noindent
6) $J \subseteq [\omega]^{\aleph_0}$ is hereditary \when \, $A \in
   [\omega]^{\aleph_0} \wedge A \subseteq^* B \in J \Rightarrow A \in J$. 

\noindent
7) $J \subseteq [\omega]^{\aleph_0}$ is dense \when \, $(\forall B \in
 [\omega]^{\aleph_0})(\exists A \in J)[A \subseteq B]$.
\end{definition}
\newpage

\section {The simple case: ${\gs} < {\ga}_*$} 

We here give a proof for the case ${\gs} < {\ga}_*$.
\begin{theorem}
\label{4d.3}  
1) If ${\gs} < {\ga}_*$ \then \, there is a saturated MAD 
family ${\cA} \subseteq [\omega]^{\aleph_0}$.

\noindent
2) Moreover, given a dense $J_* \subseteq
   [\omega]^{\aleph_0}$ we can demand ${\cA} \subseteq J_*$.
\end{theorem}

\begin{PROOF}{\ref{4d.3}}
\underline{Stage A}:  Let $\kappa = {\frak s}$, so cf$(\kappa)
> \aleph_0$.  For part (1) let $J_* \subseteq [\omega]^{\aleph_0}$ be a
dense (and even hereditary) subset of $[\omega]^{\aleph_0}$, 
i.e. as in part (2) and in both cases \wilog \, every finite union of
members of $J_*$ is co-infinite, i.e. $\omega \notin \text{ id}_{J_*}$.

Choose a sequence $\langle C^*_\alpha:\alpha < \kappa\rangle$ of
subsets of $\omega$ exemplifying ${\frak s} = \kappa$,
i.e. $\neg(\exists B \in 
[\omega]^{\aleph_0}) \wedge \bigwedge\limits_{\alpha} (B
\subseteq^* C^*_\alpha \vee B \subseteq^* \omega \backslash
C^*_\alpha)$.  For $i < \kappa$ and $\eta \in {}^i 2$ let $C^*_\eta =
C^*_i$, the aim of this notation is to simplify later proofs where we
say ``repeat the present proof but ...".
\bigskip

\noindent
\underline{Stage B}:  For $\alpha \le 2^{\aleph_0}$ let AP$_\alpha$, the set
of $\alpha$-approximations, be the set of $t$ consisting of the
following objects satisfying the following conditions:
\mn
\begin{enumerate}
\item[$\boxplus_1$]   $(a) \quad {\cT} = {\cT}_t$ is a subtree
of ${}^{\kappa >} 2$, i.e. closed under initial segments
\sn
\item[${{}}$]   $(b) \quad$ let suc$({\cT}) = \{\eta\in {\cT}:\ell
g(\eta)$ is a successor ordinal$\}$ and\footnote{so $c \ell(\{<>\}) =
\{<>,<0>,<1>\}$} 

\hskip25pt $c \ell({\cT}) = \{\eta
\in {}^{\kappa \ge}2$: if $i < \ell g(\eta)$ then $\eta \rest i \in
{\cT}\}$
\sn
\item[${{}}$]  $(c) \quad 1 \le |{\cT}| \le \aleph_0 + |\alpha|$
\sn
\item[${{}}$]   $(d) \quad \bar I = \bar I_t = \langle I_\eta:\eta \in
c \ell(\cT)\rangle = \langle I^t_\eta:\eta \in c  \ell(\cT_t)\rangle$
\sn
\item[${{}}$]  $(e) \quad \bar A = \bar A_t = \langle A_\eta:\eta
\in \text{ suc}({\cT})\rangle = \langle A^t_\eta:\eta \in \text{
suc}({\cT}_t)\rangle$ 
\end{enumerate}
\mn
such that
\mn
\begin{enumerate}
\item[${{}}$]  $(f) \quad A_\eta \in I_\eta \cap J_*$ 
or\footnote{the case ``$A_\eta = \emptyset$" is not needed in this
proof} $A_\eta = \emptyset$ and ${\mathscr S}_t = \{\eta \in 
\text{ suc}({\cT}_t):A_\eta \ne \emptyset\}$
\sn
\item[${{}}$]   $(g) \quad I_\eta = \{A \in [\omega]^{\aleph_0}$: if
$i < \ell g(\eta)$ then $A \subseteq^* 
(C^*_{\eta \rest i})^{[\eta(i)]}$ and if $i+1 < \ell g(\eta)$ 

\hskip25pt  then $A \cap A_{\eta \rest (i+1)}$ is finite$\}$,
so $I_\eta$ is well defined also when

\hskip25pt  $\eta \in c \ell({\cT})$.
\end{enumerate}
\mn
We let 
\mn
\begin{enumerate}
\item[${{}}$]   $(h) \quad C^t_\eta = C^*_\eta$ (for generalizations)
\sn
\item[$\boxplus_2$]  AP $= \cup\{\text{AP}_\alpha:\alpha \le
2^{\aleph_0}\}$
\sn
\item[$\boxplus_3$]   $s \le_{\text{AP}} t$ \Iff \, (both
are from AP and)
\begin{enumerate}
\item[$(a)$]  ${\cT}_s \subseteq {\cT}_t$
\sn
\item[$(b)$]   $\bar I_s = \bar I_t \rest c \ell(\cT_s)$
\sn
\item[$(c)$]  $\bar A_s = \bar A_t \rest \text{ suc}({\cT}_s)$.
\end{enumerate}
\end{enumerate}
\bigskip

\noindent
\underline{Stage C}:  We assert various properties of AP; of course $s,t$
denote members of AP:
\mn
\begin{enumerate}
\item[$\boxplus_4$]   $(a) \quad \le_{\text{AP}}$ partially orders AP
\sn
\item[${{}}$]   $(b) \quad \eta \triangleleft \nu \in c \ell({\cT}_t)
\Rightarrow I^t_\nu \subseteq I^t_\eta$
\sn
\item[${{}}$]   $(c) \quad$ if $\eta \in c \ell(\cT_t)$ then $I^t_\eta
\in \text{ OB}$, i.e. $I^t_\eta \cup [\omega]^{< \aleph_0}$
is an ideal of ${\cP}(\omega)$ 
\sn
\item[${{}}$]  $(d) \quad \langle A^t_\eta:\eta \in 
{\mathscr S}_t\rangle$ is almost 
disjoint (so $A^t_\eta \in \text{ ob}(\omega)$ and

\hskip25pt $\eta \ne \nu \in \mathscr{S}_t \Rightarrow A^t_\eta \cap
A^t_\nu$ finite; recall that here we can assume

\hskip25pt  ${\mathscr S}_t = \text{ suc}({\cT}_t)$)
\sn
\item[${{}}$]  $(e) \quad$ if $\eta \in c \ell({\cT}_t)$
 and $\ell g(\eta) = \kappa$ then $I^t_\eta = \emptyset$
\sn
\item[${{}}$]  $(f) \quad$ if $s \le_{\text{AP}} t$ \then \, 
$c \ell(\cT_s) \subseteq c \ell(\cT_t)$ 
and $\eta \in c \ell(\cT_s) \Rightarrow I^s_\eta = I^t_\eta$ 

\hskip25pt (and clause (b) of $\boxplus_3$ follow from clauses
(a),(c))
\sn
\item[${{}}$]  $(g) \quad \bullet \quad$ if 
$\nu \in c \ell(\cT_s) \backslash \cT_s$
and $\eta \in \mathscr{S}_s$ and $B \in I^s_\nu$ \then \, $B \cap
A_\eta$ is finite
\sn
\item[${{}}$]  $\quad \quad \bullet \quad$ if $\nu \in \cT_s$ and 
$\eta \in \mathscr{S}_s$ but $\neg(\nu \trianglelefteq \eta)$ and 
$B \in I^s_\nu$

\hskip25pt  \then \, $B \cap A_\eta$ is finite.
\end{enumerate}
\mn
[Why clause (d)?  Let $\eta_0 \ne \eta_1 \in {\mathscr S}_t$, if $\eta_0 \perp
\eta_1$ let $\rho = \eta_0 \cap \eta_1$ hence for some $\ell \in
\{0,1\}$ we have $\rho \char 94
\langle \ell \rangle \trianglelefteq \eta_0,\rho \char 94 \langle
1-\ell\rangle \trianglelefteq \eta_1$ so $A_{\eta_k} \in
I^t_{\eta_k} \subseteq I^t_{\rho \char 94 <k>} \subseteq 
\text{ ob}((C^t_\rho)^{[k]})$ for $k=0,1$ hence $A_{\eta_0} \cap A_{\eta_1}
\subseteq^* \text{ ob}((C^t_\rho)^{[\ell]}) \cap 
\text{ ob}((C^t_\rho)^{[1-\ell]}) = \emptyset$.  
If $\eta_0 \triangleleft \eta_1$ note that
$A^t_{\eta_1} \in I^t_{\eta_1} \subseteq \text{ ob}(\omega \backslash
A^t_{\eta_0})$ by clause $\boxplus_1(g)$.  
Also if $\eta_1 \triangleleft \eta_0$ similarly so
clause (d) holds indeed.

Why Clause (e)?  Recall the choice of $\langle C^*_\alpha:\alpha <
\kappa\rangle$ and $\langle C^t_\eta:\eta \in {}^{\kappa >}2 \rangle$
hence $\alpha < \kappa \Rightarrow C^t_{\eta \rest \alpha} = C^*_\alpha$.  
So if $B \in I^t_\eta$, then $B \in I_{\eta \rest(\alpha +1)}$ hence
$(B \subseteq^* C^*_\alpha \vee B 
\subseteq^* \omega \backslash C^*_\alpha)$ for every
$\alpha < \kappa$, a contradiction to the choice of $\langle
C^*_\alpha:\alpha < \kappa\rangle$.]
\mn
\begin{enumerate}
\item[$\boxplus_5$]  $(a) \quad \alpha < \beta \le 2^{\aleph_0}
\Rightarrow \text{ AP}_\alpha \subseteq \text{ AP}_\beta$
\sn
\item[${{}}$]  $(b) \quad \text{ AP}_0 \ne \emptyset$ (e.g. use $t$
with ${\cT}_t = \{<>\}$)
\sn
\item[${{}}$]  $(c) \quad$ if $\langle t_i:i < \delta\rangle$ is
$\le_{\text{AP}}$-increasing, $t_i \in \text{ AP}_{\alpha_i}$ for $i <
\delta,\langle \alpha_i:i < \delta\rangle$ is

\hskip25pt  increasing, $\delta$ a
limit ordinal and $\alpha_\delta = \cup\{\alpha_i:i < \delta\}$
\then \, 

\hskip25pt $t_\delta = \cup\{t_i:i < \delta\}$ naturally defined belongs to
AP$_{\alpha_\delta}$ and 

\hskip25pt $i < \delta \Rightarrow t_i \le_{\text{AP}} t_\delta$
\sn
\item[$\boxplus_6$]   let $J_t$ be the ideal on ${\cP}(\omega)$
generated by $\{A^t_\eta:\eta \in {\mathscr S}_t\} \cup [\omega]^{< \aleph_0}$.
\end{enumerate}
\mn
For $s \in \text{ AP}$ and $B \in \text{ ob}(\omega)$ we define:
\mn
\begin{enumerate}
\item[$(*)_1$]   $S_B = S^s_B := S^1_B \cup S^2_B$ where
\begin{enumerate}
\item[$(a)$]  $S^1_B = S^{s,1}_B :=
\{\eta \in c \ell({\cT}_s):
[B \backslash A]^{\aleph_0} \cap I^s_\eta \ne \emptyset$ for every
$A \in J_s\}$
\sn
\item[$(b)$]  $S^2_B = S^{s,2}_B := \{\eta \in c \ell(\cT_s)$: for
infinitely many $\nu,\eta \trianglelefteq \nu \in \mathscr{S}_s$ and the
set $B \cap A_\nu$ is infinite$\}$
\sn
\item[$(c)$]  $S^3_B = S^{3,s}_B := S_B$
\end{enumerate}
\sn
\item[$(*)_2$]   SP$^\iota_B = \text{ SP}^{s,\iota}_B := 
\{\eta \in {\cT}_s:\eta \char 94 \langle 0 \rangle \in S^{s,\iota}_B$ 
and $\eta \char 94 \langle 1 \rangle \in S^{s,\iota}_B\}$ for $\iota=1,2,3$
and $S_B = S^s_B = S^3_B$.
\end{enumerate}
\mn
Note
\mn
\begin{enumerate}
\item[$(*)_3$]   for $\iota=1,2,3$
\begin{enumerate}
\item[$(a)$]  $S^\iota_B$ is a subtree of $c \ell({\cT}_s)$
\sn
\item[$(b)$]  $\langle \rangle \in S_B \Leftrightarrow B \in J^+_s
\Leftrightarrow \langle \rangle \in S^1_B$
\sn
\item[$(c)$]   SP$^\iota_B \subseteq {\cT}_s$
\sn
\item[$(d)$]   if $B \subseteq A$ are from
$[\omega]^{\aleph_0}$ then $S^\iota_B \subseteq S^\iota_A$, 
SP$^\iota_B \subseteq \text{ SP}^\iota_A$. 
\end{enumerate}
\end{enumerate}
\mn
[Why?  For $\iota=1$, the first statement holds by recalling 
$\boxplus_4(b)$, the second, $\langle \rangle \in S_B
\Leftrightarrow B \in J^+_s$, holds as $I^s_{<>} = \text{ ob}(\omega)$, the
third, SP$^\iota_B \subseteq \cT_s$ as by the definition of $c \ell(\cT_s)$
we have $\eta \char 94 \langle \ell \rangle \in c \ell(\cT_s) \Rightarrow \eta
\in \cT_s$.  Also the fourth is obvious.  For $\iota=2$ this is even
easier and for $\iota=3$ it follows.]
\mn
\begin{enumerate}
\item[$(*)_4$]  if $\eta \in S_B$ and $\nu_0 \triangleleft \nu_1
\triangleleft \ldots \triangleleft \nu_{n-1}$ list $\{\nu
\triangleleft \eta:\nu \in \text{ SP}_B\}$ so this set is finite and we let
$C_s(\eta,B) := \cap\{(C^s_{\nu_\ell})^{[\eta(\ell
g(\nu_\ell))]}:\ell <n\}$, \then \, 
$S_{B \cap C_s(\eta,B)} = \{\nu \in S_B:\nu \trianglelefteq
\eta$ or $\eta \trianglelefteq \nu\}$.
\end{enumerate}
\mn
[Why?  Clearly $(\forall A \in I^s_\eta)(A \subseteq^* C_s(\eta,B))$ by
the definition of $I^s_\eta$, see $\boxplus_1(g)$ but $(\exists A
\subseteq I^s_\eta)(|A \cap B| = \aleph_0)$ hence $B \cap C_s(\eta,B) \in
\text{ ob}(\omega)$.

As $B \cap C_s(\eta,B) \subseteq B$ clearly $S_{B \cap C_s(\eta,B)}
\subseteq S_B$.  Also as $\eta \in S_B$ and as $(\forall A \in
I^s_\eta)(A \subseteq^* C_s(\eta,B))$ clearly $\eta \in S_{B \cap
C_s(\eta,B)}$ and moreover $\{\nu \in S_{B \cap C_s(\eta,B)}:\eta
\trianglelefteq \nu\} = \{\nu \in S_B:\eta \trianglelefteq \nu\}$ by
$\boxplus_4(b)$. 

Also as $S_B$ and $S_{B  \cap C_s(\eta,B)}$ are subtrees clearly
$\{\nu:\nu \trianglelefteq \eta\} \subseteq S_B \cap S_{B \cap
C_s(\eta,B)}$ and $\eta \trianglelefteq \nu \in c \ell(\cT_s)
\Rightarrow [\nu \in S_B \Leftrightarrow \nu \in S_{B \cap C_s(\eta,B)}]$.

So to prove the equality it suffices to assume $\alpha < \ell
g(\eta),\nu \in S_B,\ell g(\eta \cap \nu) = \alpha,\ell g(\nu) >
\alpha$ and $\nu \in S_{B \cap C_s(\eta,B)}$ and get a contradiction.
If $\ell < n$ and $\alpha = \ell g(\nu_\ell)$ then $(\forall A \in
I^s_\nu)[A \subseteq^* (C^s_{\eta \rest \alpha})^{[1-\eta(\alpha)]}]$, 
so an easy contradiction.  If $\alpha \notin \{\ell g(\nu_\ell):\ell < n\}$ we
can get contradiction to $\eta \rest \alpha \notin \text{\rm SP}_B$.  So we are
done proving $(*)_4$.]
\mn
\begin{enumerate}
\item[$(*)_5$]   $(a) \quad$ for every $\eta \in c \ell({\cT}_t)$
the set $\{B \in \text{ ob}(\omega):\eta \notin S_B\}$ belongs to OB
\sn
\item[${{}}$]  $(b) \quad$ if $\iota = 1,2,3$ and
$B = B_0 \cup \ldots \cup B_n
\subseteq \omega$ then $S^\iota_B = S^\iota_{B_0} \cup \ldots \cup 
S^\iota_{B_n}$
\sn
\item[${{}}$]  $(c) \quad$ if $A \in J_t$ and $B_2 = B_1 \backslash A$ then 
$S^\iota_{B_2} = S^\iota_{B_1}$, SP$^\iota_{B_2} = \text{
SP}^\iota_{B_1}$ for $\iota = 1,2,3$
\sn
\item[${{}}$]  $(d) \quad S^2_B \subseteq \cT_t$ for $B \in \text{ob}(\omega)$ 
\sn
\item[${{}}$]  $(e) \quad$ if $\eta \triangleleft \nu \in c
\ell(\cT_t)$ then $\eta \in \cT_t$
\sn
\item[${{}}$]  $(f) \quad$ if $B \in \text{ob}(\omega)$ and $s
\le_{\text{AP}} t$ \then
\begin{enumerate}
\item[${{}}$]  $\bullet_1 \quad S^{s,1}_B \supseteq S^{t,1}_B \cap c
\ell(\cT_s)$
\sn
\item[${{}}$]  $\bullet_2 \quad S^{s,2}_B \subseteq S^{t,2}_B \cap
\cT_s$ (inclusion in different direction?  yes!)
\sn
\item[${{}}$]  $\bullet_3 \quad S^{s,3}_B \supseteq S^{t,3}_B \cap c
\ell(\cT_s)$, in fact $(S^{t,2}_B \cap c \ell(\cT_s) \backslash
S^{s,2}_B \subseteq S^{s,1}_B$
\sn
\item[${{}}$]  $\bullet_4 \quad$ if $\eta \in S^{s,1}_B \backslash
\cT_t$ then $\eta \in S^{t,1}_B$

\hskip25pt  (clause (f) is not used here)
\end{enumerate}
\item[$(*)_6$]   for $\iota = 1,2,3$ we have
\begin{enumerate}
\item[$(a)$]   if $B \subseteq \omega,\ell < 2$ 
and $\nu \in S^\iota_B \cap {\cT}_t$ and $B \subseteq^* (C^t_\nu)^{[\ell]}$
\then \, $\nu \char 94 \langle \ell \rangle \in S^\iota_B$
but $\nu \char 94 \langle 1-\ell \rangle \notin S^\iota_B$
\sn
\item[$(b)$]   if $B \subseteq \omega$ 
and $\nu \in S^\iota_B \cap {\cT}_t$ \then \, for some $\ell < 2$ we have 
$\nu \char 94 \langle \ell \rangle \in S^\iota_B$
\sn
\item[$(c)$]   $\omega \notin J_t$.
\end{enumerate}
\end{enumerate}
\mn
[Why?  Read the definitions recalling $(*)_5(c)$.  For clause (c)
recall that $\nu \in \mathscr{S}_s \Rightarrow A_s \in J_s$ and by
$\boxplus_1(f)$ we have $J_s \subseteq J_*$ and by Stage A, $\omega
\notin \text{ ob}(J_*)$.]
\bigskip

\noindent
\underline{Stage D}:
\mn
\begin{enumerate}
\item[$\boxplus_7$]   if $\alpha < 2^{\aleph_0},s \in \text{
AP}_\alpha$ and $B \in \text{ ob}(\omega) \backslash J_s$ \then \,
we can find $t \in \text{ AP}_{\alpha +1}$ such that $s
\le_{\text{AP}} t$ and $B$ contains $A_\eta$ for some $\eta \in 
{\mathscr S}_t \backslash {\cT}_s$.
\end{enumerate}
\mn
This is a major point and we shall prove it in Stage F below.
\bigskip

\noindent
\underline{Stage E}:  We prove the theorem.

Let $\langle B_\alpha:\alpha < 2^{\aleph_0}\rangle$ list ${\cP}(\omega)$
each appear $2^{\aleph_0}$ times.  By induction on $\alpha \le
2^{\aleph_0}$ we choose $t_\alpha$ such that
\mn
\begin{enumerate}
\item[$\circledast$]  $(a) \quad t_\alpha \in \text{ AP}_\alpha$
\sn
\item[${{}}$]  $(b) \quad \beta < \alpha \Rightarrow t_\beta
\le_{\text{AP}} t_\alpha$
\sn
\item[${{}}$]  $(c) \quad$ if $\alpha = \beta +1$ then either
$B_\alpha \in J_{t_\beta}$ or $B_\alpha$ contains $A_\eta$, for 

\hskip25pt some $\eta \in {\mathscr S}_{t_\alpha} \backslash {\cT}_{t_\beta}$.
\end{enumerate}
\mn
For $\alpha = 0$ use $\boxplus_5(b)$.

For $\alpha$ limit use $\boxplus_5(c)$.

For $\alpha = \beta +1$ use $\boxplus_7$.

Now let $t \in \text{ AP}$ be $\cup\{t_\alpha:\alpha < 2^{\aleph_0}\}$
and recalling $(*)_6(c)$ it is easy to check that $\bar A_t$ is a 
saturated MAD family,
enough for \ref{4d.3}(1) and recalling that by $\boxplus_1(f)$ it is
$\subseteq J_*$ also enough for \ref{4d.3}(2).
\bigskip

\noindent
\underline{Stage F}:  The rest of the proof is dedicated to

\underline{the proof of $\boxplus_7$ so $\alpha,s$ and $B$ are given}.

The proof is now split into cases. 
\bigskip

\noindent
\underline{Case 1}:  Some $\nu \in S_B$ is such that 
$\nu \in c \ell({\cT}_s) \backslash {\cT}_s$.

By $(*)_5(d)$ we have $\nu \in S^1_B$.
Clearly as $\nu \in S_B$ there is $B_1 \in [B]^{\aleph_0} \cap
I^s_\nu$.  Note that $\ell g(\nu) > 0$ as $\langle \rangle \in
{\cT}_s$ by clause (c) of $\boxplus_1$.

Note that $A \in I^s_\nu \wedge \eta \in {\mathscr S}_s \Rightarrow A \cap
A^s_\eta$ is finite, e.g. by the proof of $\boxplus_4(d)$ or better by
$\boxplus_4(g)$.
\bigskip

\noindent
\underline{Subcase 1A}:  Assume $\ell g(\nu)$ is a successor ordinal.  

Let $B_2 \subseteq B_1$ be such that $B_2 \in J_*$ and $B_1 \backslash B_2$ are
infinite.  Now define $t$ as follows: ${\cT}_t = {\cT}_s \cup
\{\nu\},A^t_\rho$ is $A^s_\rho$ if $\rho \in \text{ suc}({\cT}_s)$
and is $B_2$ if $\rho = \nu$, lastly define $I^t_\rho$ for $\rho \in
{\cT}_t$ as in clause (g) of $\boxplus_1$.  Easy to check 
that $t$ is as required; actually $B_2 = B_1$ is O.K., too.
\bigskip

\noindent
\underline{Subcase 1B}:  Assume $\ell g(\nu)$ is a limit ordinal.

Clearly $\ell g(\nu) < \kappa$ by $\boxplus_4(e)$, as $I^s_\nu \ne
\emptyset$ because $B_1 \in I^s_\nu$, clearly 
there is $\ell \in \{0,1\}$ such that $B'_1 :=
(C^s_\nu)^{[\ell]} \cap B_1$ is infinite, let $B_2 \subseteq
B'_1$ be such that $B_2,B'_1 \backslash B_2$ are infinite and $B_2 \in
J_*$.  We define
$t$ by ${\cT}_t = {\cT}_s \cup \{\nu,\nu \char 94 \langle \ell
\rangle\},A^t_\rho$ is $A^s_\rho$ if $\rho \in \text{ suc}({\cT}_s)$ 
and is $B_2$ if $\rho = \nu \char 94 \langle \ell \rangle$ and
$I^t_\rho$ for $\rho \in {\cT}_t$ is defined 
as in clause (g) of $\boxplus_1$.

Easy to check that $t$ is as required.
\bigskip

\noindent
\underline{Case 2}:  SP$_B = \emptyset$ but not case 1.

Let $\nu^*_B := \cup\{\eta:\eta \in S_B\}$.
\bigskip

\noindent
\underline{Subcase 2A}:  $\nu^*_B \in S^1_B$.

As $S_B \subseteq c \ell({\cT}_s)$ by the definition of $S_B$,
as we are assuming ``not case 1" necessarily $S_B \subseteq \cT_s$ hence
$\nu^*_B \in {\cT}_s$ so $\ell g(\nu^*_B) < \kappa$.  

We define $B^*_2$ as $B \cap A_{\nu^*_B}$ if $A_{\nu^*_B}$ is well
defined and $B^*_2 = \emptyset$ otherwise; and for $\ell=0,1$, 
let $B^*_\ell := B \cap (C^*_{\ell g(\nu^*_b)})^{[\ell]} \backslash B^*_2$.

So
\mn
\begin{enumerate}
\item[$(*)_7$]   $\langle B^*_0,B^*_1,B^*_2 \rangle$ is a partition of $B$ 
\end{enumerate}
\mn
hence by $(*)_5(b)$ for some $\ell=0,1,2$
\mn
\begin{enumerate}
\item[$(*)_8$]   $\nu^*_B \in S^1_{B^*_\ell}$ 
\end{enumerate}
\mn
easily
\mn
\begin{enumerate}
\item[$(*)_9$]   $\ell \ne 2$,
\end{enumerate}
\mn
and
\mn
\begin{enumerate}
\item[$(*)_{10}$]   $\rho := \nu^*_B \char 94 \langle \ell
\rangle \in S_{B^*_\ell}$.
\end{enumerate}
\mn
[Why? By the definitions noting $\rho \in c \ell(\cT_s)$.]

Also as $B^*_\ell \subseteq B$ clearly
\mn
\begin{enumerate}
\item[$(*)_{11}$]   $S_{B^*_\ell} \subseteq S_B$.
\end{enumerate}
\mn
But $(*)_{10} + (*)_{11}$ contradicts the choice of $\nu^*_B$.
\bigskip

\noindent
\underline{Subcase 2B}:  $\nu^*_B \notin S^1_B$.

By $(*)_3(b) + (*)_6(c)$ 
and the assumption of $\boxplus_7$ we have $\langle \rangle
\in S_B$ and by $(*)_6(b)$ 
clearly $\langle 0 \rangle \in S_B$ or $\langle 1 \rangle
\in S_B$ hence $\nu^*_B \ne \langle \rangle$.  If $\nu^*_B = \nu \char
94 \langle \ell \rangle$ by the definition of $\nu^*_B$ 
we have $\nu^*_B \in S_B$, contradiction to the subcase assumption.  
Hence necessarily $\ell g(\nu^*_B)$
is a limit ordinal $\le \kappa$, call it $\delta$.  So $\alpha <
\delta \Rightarrow \nu^*_B \rest \alpha \in S_B$ but $\rho
\triangleleft \varrho \in c \ell({\cT}_s) \Rightarrow \rho \in
{\cT}_s$ hence $\alpha < \delta \Rightarrow \nu^*_B \rest \alpha \in
{\cT}_s$.  Now for every
$\alpha < \delta$ let $\nu^*_{B,\alpha} := (\nu^*_B \rest \alpha)
\char 94 \langle 1 - \nu^*_B(\alpha)\rangle$, so clearly
$\nu^*_{B,\alpha} \in c \ell({\cT}_s) \backslash S_B$ hence as
$\nu^*_{B,\alpha} \notin S^2_B$, by the definition of $S^2_B
\subseteq S_B$ the set $\cA_\alpha = \{\rho \in
\mathscr{S}_s:\nu^*_{B,\alpha} \trianglelefteq \rho$ and $B \cap
A_\rho$ is infinite$\}$ is finite, so we can
find $n=n(\alpha) < \omega$ and $A^*_{\alpha,0},
\dotsc,A^*_{\alpha,n(\alpha)-1}$ enumerating $\cA_\alpha$ but
also $\nu^*_{B,\alpha} \notin S^1_B \subseteq S_B$, hence ob$(B \backslash \cup
\cA_\alpha) = \text{ ob}(B \backslash \cup\{A^*_{\alpha,\ell}:\ell <
n(\alpha)\}$ is disjoint to $I_{\nu^*_{B,\alpha}} \cap J^+_s$ and by
the choice of $\cA_\alpha$ and $(*)_4$, ob$(B \backslash \cup \cA_\alpha)=
[B \backslash (A^*_{\alpha,0} \cup
\ldots \cup A^*_{\alpha,n(\alpha)-1})]^{\aleph_0}$ is disjoint to
$I^s_{\nu^*_{B,\alpha}}$.
Let $A^*_{\alpha,n(\alpha)}$ be
$A_{\nu^*_B \rest \alpha}$ when defined and $\emptyset$ otherwise.  By the
definitions of $I^s_{\nu^*_{B,\alpha}},I^s_{\nu^*_B \rest \alpha}$ we
have (for $\alpha < \delta$ of course):
\mn
\begin{enumerate}
\item[$\odot_1$]   $(a) \quad [B \cap (C^s_{\nu^*_B \rest \alpha})
^{[1-\nu^*_B(\alpha)]} \backslash (A^*_{\alpha,0} \cup
\ldots \cup A^*_{\alpha,n(\alpha)})]^{\aleph_0}$ is disjoint to
$I^s_{\nu^*_B \rest \alpha}$
\sn
\item[${{}}$]  $(b) \quad A^*_{\alpha,\ell} \in \{A_\rho:\rho \in
\mathscr{S}_t$ and $\nu^*_{B,\alpha} \trianglelefteq \rho$ (hence
$A_{\alpha,\ell} \subseteq (C^s_{\nu^*_B \rest 
\alpha})^{[1-\nu^*_B(\alpha)]})\}$ 

\hskip25pt for $\ell < n(\alpha)$ (not needed presently)
\end{enumerate}
\mn
[Why clause $(b)$?   By the choice of $\cA_\alpha$.]

Let ${\cA}^* = \{B \cap A:A = A^*_{\alpha,k}$ for some $\alpha <
\delta,k \le n(\alpha)$ and $B \cap A$ is infinite$\}$.

So ${\cA}^*$ is a family of pairwise almost disjoint infinite
subsets of $B$ and if ${\cA}^*$ is finite, still $B \backslash
\cup\{A:A \in {\cA}^*\}$ is infinite because ${\cA}^* \subseteq
J_s$ and we are assuming $B \notin J_s$.  

Let $\Lambda := \{\nu \in \mathscr{S}_s:
\nu^*_B \trianglelefteq \nu,|A_\nu \cap B| = {\aleph_0}\}$.

Now
\mn
\begin{enumerate}
\item[$\odot_2$]  there is  a set $ B_ 1 $ such that: 
\begin{enumerate}
\item[$(a)$]  $B_1 \subseteq B$ is infinite 
\sn
\item[$(b)$]  $B_1$ is almost disjoint to any $A \in {\cA}^*$
\sn
\item[$(c)$]  if $\Lambda$ is finite then $\nu \in \Lambda 
\Rightarrow |B_1 \cap  A_\nu| < {\aleph_0}$ 
\sn
\item[$(d)$]  if $\Lambda$ is infinite then for infinitely many 
$\nu \in \Lambda$ we have $|B_1 \cap A_\nu| = {\aleph_0}$ 
\end{enumerate}
\end{enumerate}
\mn
[Why?  First assume $ \Lambda $ is finite, so \wilog \, it is empty.  
If ${\cA}^*$ is finite use the paragraph above on $\cA^*$.  Otherwise as
$|{\cA}^*| \le |\delta| + \aleph_0 \le \kappa = {\frak s}$ and by 
the theorem's assumption ${\frak s} < {\frak a}_* \le {\frak a}$
and by the definition of ${\frak a}$ it follows that $\odot_2$ holds. 

Second, assume that $\Lambda$ is infinite, and choose pairwise 
distinct $\nu_n \in \Lambda$ for $n < \omega$.  Now we recall that we are 
assuming ${\frak s} < {\frak a}_*$ and apply the definition 
\ref{0z.3d} of ${\frak a}_*$ to ${\cA}^*$ and 
$\langle \nu_n:n < \omega \rangle$ and get an infinite  
$B_1 \subseteq B$  as required.]
\mn
\begin{enumerate}
\item[$\odot_3$]   $(a) \quad B_1 \in J^\perp_s$
\sn
\item[{{}}]   $(b) \quad$ if $\neg(\nu^*_B \trianglelefteq \eta)$ 
and $\nu \in \mathscr{S}_s$ then $B_1 \cap A_\nu$ is finite.
\end{enumerate}
\mn
[Why? For clause (a), note that first $B_1 \subseteq B \subseteq
\omega$, second $B_1$ is infinite by clause (a) of $\odot_2$,
third $B_s \notin J_s$ is proved by dividing to two cases.  If 
$\Lambda$ is finite use clause (b) of $\odot_3$ proved below and
clause (b) of $\odot_2$; and if 
$\Lambda$ is infinite use $\odot_2(d)$). So let us turn to proving 
clause (b); we should prove that $[\eta \in \text{ suc}(\cT_s) \wedge 
\neg(\nu^*_B \trianglelefteq \eta) \Rightarrow B_1 \cap A_\eta 
\text{ finite}$.

If $A_\eta \in \{A^*_{\alpha,n}:\alpha < \delta,n \le n(\alpha)\}$ then either
$A_\eta \cap B$ is finite hence $A_\eta \cap B_1 \subseteq A_\eta \cap B$
is finite or $A_\eta \cap B$ is infinite hence $A_\eta \cap B \in
\cA^*$ hence $B_1 \cap (A_\eta \cap B)$ is finite by the choice of
$B_1$ but $B_1 \subseteq B$ hence $B_1 \cap A_\eta$ is finite.
So assume $A_\eta \notin \{A^*_{\alpha,m}:\alpha < \delta,n \le
n(\alpha)\}$, so by the choice of $A_{\alpha,n(\alpha)}$ for $\alpha <
\delta$ necessarily $\neg(\eta \triangleleft \nu^*_B)$.  
Recall that we are assuming that $\neg(\nu^*_B \trianglelefteq \eta)$.  
Together for some $\alpha < \delta$ we have
$\alpha = \ell g(\nu^*_B \cap \eta) < \delta$ and $\nu^*_B \rest
\alpha \triangleleft \eta$ and we get contradiction by the choice of
$\cA_\alpha = \{A^*_{\alpha,\ell}:\ell < n(\alpha)\}$ and
$A^*_{\alpha,n(\alpha)}$.]

We shall now prove by induction on $\alpha \le \delta$ that $B_1 \in
I^s_{\nu^*_B \rest \alpha}$.  For $\alpha=0$ recall $I^s_{\nu^*_B \rest
\alpha} = [\omega]^{\aleph_0}$, for $\alpha$ limit
$I^s_{\nu^*_B \rest \alpha} = 
\cap\{I^s_{\nu^*_B \rest \beta}:\beta < \alpha\}$ and use the
induction hypothesis.  For $\alpha = \beta +1$ first note that $B_1$
is almost disjoint to $A_{\nu^*_B \rest \beta}$ if $\nu^*_B \rest \beta
\in \mathscr{S}_s \subseteq \text{ suc}(\cT_s)$ by $\odot_2(b)$ 
and, second, $B_1$ is almost
disjoint to $(C^s_{\nu^*_B \rest \beta})^{[1-\nu^*_B(\beta)]}$ otherwise
recalling $\odot_3(b)$ we get contradiction to the 
present case assumption SP$_B = \emptyset$ by $(*)_6(a) \,+$ the induction
hypothesis; together by the definition of $I^s_{\nu^*_B \rest
\beta},I^s_{\nu^*_B \rest \alpha}$ we have $B_1 \in I^s_{\nu^*_B \rest
\alpha}$.  Having carried the induction, in
particular $B_1 \in I^s_{\nu^*_B \rest \delta} = I_{\nu^*_B}$; now
recalling first $B_1 \notin J_s$ by $\odot_3(a)$, second $B_1 \subseteq B$ by
$\odot_2(a)$ and third, the choice 
of ${\cA}^*,J_s$, together they contradict the subcase 
assumption $\nu^*_B \notin S^1_B$.
\bigskip

\noindent
\underline{Case 3}:  None of the above.

Without loss of generality
\mn
\begin{enumerate}
\item[$\oplus_1$]  if $B_1 \subseteq B$ but $B_1 \notin J_s$ then
none of the two cases above holds.
\end{enumerate}
\mn
We try to choose $\bar\eta^n = \langle \eta_\rho:\rho \in
{}^n2\rangle$ by induction on $n$ such that:
\mn
\begin{enumerate}
\item[$(a)$]  $\eta_\rho \in \text{ SP}_B$
\sn
\item[$(b)$]   if $\rho = \varrho \char 94 \langle \ell \rangle$
then $\eta_\varrho \char 94 \langle \ell\rangle \trianglelefteq
\eta_\rho$
\sn
\item[$(c)$]   $\{\nu:\nu \triangleleft \eta_\rho$ and $\nu \in
\text{ SP}_B\} = \{\eta_{\rho \rest k}:k < \ell g(\rho)\}$.
\end{enumerate}
\mn
For $n=0$, note that SP$_B \ne \emptyset$ as not case 2 (and not 
case 1) so we can
choose $\eta_\rho \in SP_B$ with minimal length.  If $n=m+1$ and $\rho
\in {}^m 2$ by the induction hypothesis $\eta_\rho \in \text{ SP}_B$, hence
$\eta_\rho \in {\cT}_s$ and by the definition of SP$_B$ for
$\ell=0,1$ the sequence $\eta_\rho \char 94 \langle \ell \rangle$
belongs to $S_B$.  

First assume $\{\nu \in \text{ SP}_B:
\eta_\rho \char 94 \langle \ell \rangle \trianglelefteq \nu\} = \emptyset$.  So
$B_1 := B \cap C_s(\eta_\rho \char 94 \langle \ell \rangle,B) 
\notin J_s$ noting $C_s(\eta_\rho \char 94 \langle \ell \rangle,B) =
\cap\{C^{[\rho(k)]}_{\eta_\rho \rest k}:k \le m\}$, recalling it
is defined in $(*)_4$ from Stage C using $\eta_\rho \char 94 \langle
\ell \rangle \in S_B$;
hence $\eta_\rho \char 94 \langle \ell \rangle \in S_{B_1}$.

Now by $(*)_4$ we know $S_{B_1} = \{\nu \in S_B:\nu \trianglelefteq
\eta_\rho \char 94 \langle \ell \rangle$ or $\eta_\rho \char 94
\langle \ell \rangle \trianglelefteq \nu\}$ so case 2 or case 1 holds
for $B_1$, contradiction to $\oplus_1$.

Second, assume we have $(\exists \eta)(\eta_\rho \char 94 \langle \ell
\rangle \trianglelefteq \eta \in \text{ SP}_B)$ so choose such $\eta_{\rho
\char 94 \langle \ell\rangle}$ of minimal length.

Hence we have carried the inductive choice of $\langle \bar \eta^n:n <
\omega\rangle$. 

For each $\rho \in {}^\omega 2$ let $\eta_\rho = \cup\{\eta_{\rho
\rest n}:n < \omega\}$, clearly $\eta_\rho \in c \ell({\cT}_s)$.
Also $\langle \eta_\rho:\rho \in {}^\omega 2\rangle$ is without
repetitions and each $\eta_\rho$ belongs to $c \ell({\cT}_s)$, so
as $|{\cT}_s| < 2^{\aleph_0}$ there is $\rho \in {}^\omega 2$ such
that $\eta_\rho \notin {\cT}_s$.  By clause (c) above we have
$\{\varrho:\varrho \triangleleft \eta_\rho$ and $\varrho \in SP_B\} =
\{\eta_{\rho \rest n}:n < \omega\}$.

Note that
\mn
\begin{enumerate}
\item[$\oplus_2$]   $\langle C_s(\eta_{\rho \rest k},B):k < \omega\rangle$ is
$\subseteq$-decreasing.  
\end{enumerate}
\mn
Let $\cW = \{\alpha < \ell g(\eta_\rho)$: for some $\nu \in 
\mathscr{S}_s$ we have 
$\ell g(\nu \cap \eta_\rho) = \alpha$ and $A_\nu \cap B$ is infinite$\}$.

First, assume $\cW$ is an unbounded subset of 
$\ell g(\eta_\rho)$.  In this case choose $\alpha_n \in \cW$ such that 
$\alpha_{n+1} > \alpha_n \ge \ell g(\eta_{\rho \rest n})$ for
$n < \omega$ and we choose $\nu_n \in \mathscr{S}_s$ such that 
$\ell g(\nu_n \cap \eta_\rho) = \alpha_n$ and $A_{\nu_n} \cap B$ is infinite.
So can choose an infinite $B_0 \subseteq B$ such that
$n < \omega$ implies $B_0 \backslash \cup\{A_{\eta_{\rho \rest k}}:
k < n\} \subseteq^* C_s(\eta_{\rho \rest \alpha_n},B)$ and 
$(B_0 \cap A_{\nu_n} \in \text{ ob}(\omega))$.

So
\mn
\begin{enumerate}
\item[$\oplus_3$]    $B_0 \subseteq B,B_0 \notin J_s$;
\sn
\item[$\oplus_4$]    the set SP$_{B_0}$ is empty.
\end{enumerate}
\mn
[Why?  By $(*)_4$ for each $n < \omega$ we have $S_{B_0} \subseteq \{\nu:\nu
\trianglelefteq \eta_{\rho \rest n} \vee \eta_{\rho \rest n} 
\trianglelefteq \nu\}$, hence
$S_{B_0} \cap {\cT}_s \subseteq \{\nu:\nu \triangleleft \eta_\rho$ or
$\eta_\rho \trianglelefteq \nu\}$ but $\eta_\rho \notin \cT_s$ 
so SP$_{B_0} = \emptyset$.]
\mn
\begin{enumerate}
\item[$\oplus_5$]    $S_{B_0}$ is not empty.
\end{enumerate}
\mn
[Why?  By $\boxplus_3$.]

By $\oplus_4 + \oplus_5$ for the set $B_0$, case 2 or case 1 hold,
 so we get contradiction to $\oplus_1$.

Second, assume $\sup(\cW ) < \ell g(\eta_\rho)$,
so we can choose $n(*) < \omega$ such that $\sup(\cW ) < 
\ell g(\eta_{\rho \rest n(*)})$.
Now $[\nu \in \mathscr{S}_s \wedge \eta_{\rho \rest n(*)}
\trianglelefteq \nu \Rightarrow B \cap A^s_\rho$ is finite] as
otherwise recalling 
$\eta_\rho \in c \ell(\cT_s) \backslash \cT_s$ necessarily
$\alpha = \ell g(\eta_\rho \cap \nu) < \ell g(\eta_\rho)$ and of
course $\alpha \ge \ell g(\eta_{\rho \rest n(*)})$, but see the choice
of $n(*)$; so $\eta_{\rho \rest n(*)} \notin S^2_\beta$ hence
$\eta_{\rho \rest n(*)} \in S^1_B$, so we can 
choose an infinite $B_1 \subseteq B$ such that
$ B_1 \in I^s_{\eta_{ \rho \rest n(*)}}$.  So checking by cases, 
$B_1 \in {\rm ob} (\omega)$ is almost disjoint to any $A_\nu,\nu \in
\mathscr{S}_s$.  Obviously $B_1 \in I^s_{\eta_\rho}$, so for it
case 1 holds as exemplified by $\eta_\rho$ 
again contradiction to $\oplus_1$.  
\end{PROOF}
\newpage

\section {The other cases} 

\begin{theorem}
\label{g.3} 
1) If $\kappa = {\gs} = {\ga}_*$ and 
{\rm cf}$([{\gs}]^{\aleph_0},\subseteq) = {\gs}$ \then \, there
   is a saturated MAD family.

\noindent
2) If $\kappa = {\gs} = {\ga}_*$ and $\bold U(\kappa) = \kappa$, see
Definition \ref{g.5} below and $J_* \subseteq [\omega]^{\aleph_0}$ is 
dense \underline{then} there is a saturated 
MAD family $\subseteq J_*$.
\end{theorem}

\noindent
Recall
\begin{definition}
\label{g.5} 
1) For cardinals $\partial \le \sigma \le \theta \le \lambda$ 
(also the case $ \theta < \sigma $ is OK)
let $\bold U_{\theta,\sigma,\partial}(\lambda) = \text{ Min}
\{|{\cP}|:{\cP} \subseteq [\lambda]^{\le \sigma}$ such that for
every $X \in [\lambda]^\theta$ for some $u \in {\cP}$ we have $|X
\cap u| \ge \partial\}$.  If $\partial = \sigma$ we may omit
$\partial$; if $\sigma = \partial = \aleph_0$ we may omit them both, and if
$\sigma = \partial = \aleph_0 \wedge \theta = \lambda$ we may omit 
$\theta,\sigma,\partial$.
In the case of our Theorem, it means: $\bold U(\kappa) = \text{
Min}\{|{\cP}|:{\cP} \subseteq [\kappa]^{\le \aleph_0}$ and $(\forall X
\in [\kappa]^\kappa)(\exists u \in {\cP})(|X \cap u| \ge \aleph_0)\}$.

\noindent
2) If in addition $J$ is an ideal on $\theta$ \then \,
$\bold U_{\theta,\sigma,J}(\lambda) = \text{ Min}\{|{\cP}|:{\cP}
\subseteq [\lambda]^{\le\sigma}$ such that for every function $f:\theta
\rightarrow \lambda$ for some $u \in {\cP}$ the set $\{i <
 \theta:f(i) \in u\}$ does not belong to $J\}$. 

\noindent
3) Let Pr$(\kappa,\theta,\sigma,\partial)$ mean: $\kappa \ge \theta
 \ge \sigma \ge \partial$ and we can find $(E,\bar{\cP})$ 
such that (if $\partial = \aleph_0$ we may omit
 $\partial$, if $\sigma = \partial = \aleph_0$ we may omit them, if
 $\sigma = \partial = \aleph_0 \wedge \theta = \kappa$ we may omit
 $\theta,\sigma,\partial$):
\mn
\begin{enumerate}
\item[$(a)$]   $\bar{\cP} = \langle \cP_\alpha:\alpha \in E \rangle$
\sn
\item[$(b)$]   $E$ is a club of $\kappa$ and $\gamma \in E \Rightarrow
|\gamma|$ divide $\gamma$
\sn
\item[$(c)$]   if $u \in \cP_\alpha$ then $u \in [\alpha]^{\le
\sigma}$ has no last member
\sn
\item[$(d)$]   $\bullet_1 \quad \bar{\cP}$ is $\subseteq$-increasing
\sn
\item[${{}}$]  $\bullet_2 \quad |\cP_\alpha| < \kappa$
\sn
\item[$(e)$]  if $w \subseteq \kappa$ is bounded and otp$(w) = \theta$ and
sup$(w) \in \text{ acc}(E)$  \then \, for some $u,j$ we have:
\sn
\item[${{}}$]  $\bullet_1 \quad |u \cap w| \ge \partial$
\sn
\item[${{}}$]  $\bullet_2 \quad j \in \text{ acc}(E)$
\sn
\item[${{}}$]  $\bullet_3 \quad u \in \cP_j$
\sn
\item[${{}}$]  $\bullet_4 \quad |w \cap j| < \theta$,i.e. $j < \sup(w)$
\sn
\item[$(f)$]   if $i \in \{0\} \cup E$ and $j = 
\text{ min}(E \backslash (i+1)),w \subseteq [i,j)$, otp$(w) = 
\theta$ \then \, for some set $u$
\sn
\item[${{}}$]  $\bullet_1 \quad u \in \cP_j$ and $u \subseteq (i,j)$
\sn
\item[${{}}$]  $\bullet_2 \quad |u \cap w| \ge \partial$.
\end{enumerate}
\end{definition}

\begin{explanation}
\label{g.7}  
The proof of \ref{g.3} is based on the proof of \ref{4d.3}.  
The difference is that in the proof
of $\odot_2$ of subcase 2B of stage F, if $\ell g(\nu^*_B) = \kappa$ it
does not follow that we have $|{\cA}^*| < {\ga}_*$, 
so we have to do something else when $|{\cA}^*| = {\ga}_* = {\gs}$.  
By the assumption $\bold U(\kappa) = \kappa$ there is a sequence 
$\langle u_\alpha:\omega \le \alpha < \kappa\rangle$
of members of $[\kappa]^{\aleph_0}$ such that $u_\alpha \subseteq
\alpha$ and for every $X \in [\kappa]^\kappa$ for some
$\alpha,u_\alpha \cap X$ is infinite.  Now if e.g.
$\ell g(\nu) =
\alpha \ge \omega$ we can use $u_\alpha$ and apply \ref{g.19} below
to appropriate $\bar B_\nu$ and get ${\cP}_\nu$ 
and add it to the family $\{C^*_\alpha:\alpha
< \kappa\}$ witnessing ${\frak s} = \kappa$ the family ${\cP}_\nu$ as in
\ref{g.19}.  So now we really need to use $C^s_\nu$ rather
than $C^*_\alpha$.
\end{explanation}

\begin{observation}
\label{g.10}
If Pr$(\kappa,\theta,\sigma,\partial)$ is satisfied by
   $(E,\bar{\cP})$ \then \, we can find $(E',\bar{\cP}')$ as in
   \ref{g.5}(3) but
\mn
\begin{enumerate}
\item[$(d)'$]  $\bullet_2 \quad$ if $j > \sup(j \cap E')$ then
$|\cP'_j| \le j$
\sn
\item[$(e)$]   as above but sup$(w) \in E$.
\end{enumerate}
\end{observation}

\begin{PROOF}{\ref{g.10}}
Use any club $E' \subseteq \text{ acc}(E)$ of
$\kappa$ such that $\delta \in E' \Rightarrow |\cP_\delta| \le
|\text{min}(E' \backslash (\delta +1))|$ and
$\delta \in \text{ nacc}(E') \Rightarrow \text{
cf}(\delta) \ne \text{ cf}(\theta)$ and let $\cP'_\gamma$ be
$\cP_\gamma$ if $\gamma \in \text{ acc}(E')$ and be
$\cup\{\cP_\beta:\beta \in E \cap \gamma\}$ if $\gamma \in \text{
nacc}(E')$.
\end{PROOF}

\begin{observation}
\label{g.19}  
Assume $\bar B^* = \langle B^*_n:n <
\omega\rangle$ satisfies $B^*_n \in [\omega]^{\aleph_0},B^*_{n+1}
\subseteq B^*_n$ and $|B^*_n \backslash B^*_{n+1}| = \aleph_0$ for
infinitely many $n$'s.  \underline{Then} we can
find ${\cP}$ such that
\mn
\begin{enumerate}
\item[$(*)$]  $(a) \quad {\cP} \subseteq [\omega]^{\aleph_0}$ is
of cardinality ${\gb}$
\sn
\item[${{}}$]  $(b) \quad$ if ${\cA} \subseteq
[\omega]^{\aleph_0}$ is an AD family, $B \subseteq \omega$ and
$(\forall n)(B \cap B^*_n$ 

\hskip25pt  $\notin \text{\rm id}_{\cA})$
\then \, for some countable (infinite) ${\cP}' \subseteq {\cP}$ 
for $2^{\aleph_0}$ function 

\hskip25pt   $\eta \in {}^{{\cP}'}2$ we have:
for some {\rm id}$_{\cA}$-positive set $A \subseteq^* B$ we have: 

\hskip25pt  $A \subseteq^* C^{[\eta(C)]}$ for every $C \in {\cP}'$ and $A
\subseteq^* B_n$ for every $n$.
\end{enumerate}
\end{observation}

\begin{PROOF}{\ref{g.19}}
\underline{Proof of \ref{g.19}}  
Let ${\cB} = \{\bar B:\bar B = \langle B_n:n <\omega\rangle$ where $B_n
\subseteq \omega$ is infinite, $B_n \supseteq B_{n+1}$ and $B_n
\backslash B_{n+1}$ is infinite for infinitely many $n < \omega\}$, 
i.e. the set of $\bar B$ satisfying the demands on $\bar B^*$.

For $\bar B \in {\cB}$ and $\cA \subseteq [\omega]^{\aleph_0}$
let pos$(\bar B,{\cA}) = \{B \subseteq \omega:B \cap B_n \notin 
\text{ id}_{\cA}$ for every $n\}$.
So the claim says that for every $\bar B \in {\cB}$
there is ${\cP} \subseteq [\omega]^{\aleph_0}$ of cardinality
${\frak b}$ such that if ${\cA} \subseteq [\omega]^{\aleph_0}$ is
an AD family and $B \in \text{ pos}(\bar B,{\cA})$ \then \, there
is a countable infinite ${\cP}' \subseteq {\cP}$ as there.

Consider the statement:
\mn
\begin{enumerate}
\item[$\boxplus$]   if $\bar B \in {\cB}$ \then \, we can
find $\bold B$ such that
\begin{enumerate}
\item[$(a)$]   $\bold B = \langle \bar B_\delta:\delta \in
S^{\gb}_{\aleph_0}\rangle$ recalling $S^{\gb}_{\aleph_0} =
\{\delta < {\gb}:\text{cf}(\delta) = \aleph_0\}$
\sn
\item[$(b)$]  $\delta \in S^{\gb}_{\aleph_0} \Rightarrow
\bar B_\delta \in {\cB}$
\sn
\item[$(c)$]  if ${\cA}$ is an AD family and $B \in \text{
pos}(\bar B,{\cA})$, then for some club $E$ of ${\gb}$, for every
$\delta \in E \cap S^{\gb}_{\aleph_0}$ we have $(\exists^\infty n)[B \cap
(B_{\delta,n} \backslash B_{\delta,n+1}) \in \text{ id}^+_{\cA}]$
\sn
\item[$(d)$]  if $\delta_1 < \delta_2$ are from $S^{\gb}_{\aleph_0}$
\then \, for some $n < \omega$ the set $B_{\delta_1,n} \cap
B_{\delta_2,n}$ is finite.
\end{enumerate}
\end{enumerate}
\mn
Why is this statement enough?  By it we can find a subset $\cB'$ of
$\cB$ of cardinality $\gb$ such that $\bar B^* \in \cB'$ and for every $\bar B
\in \cB'$ for some $\bold B = \langle \bar B_\delta:\delta \in
S^{\gb}_{\aleph_0}\rangle$ as in $\boxplus$ we have $\delta \in
S^{\gb}_{\aleph_0} \Rightarrow \bar B_\delta \in \cB'$.  
Now $\cP$, the closure by Boolean operations of 
$\{B_n:\bar B \in \cB'$ and $n < \omega\}$ is as required. 

Why?  Let $\bar B \in \cB'$ (e.g. $\bar B^*$) and an
AD family $\cA \subseteq [\omega]^{\aleph_0}$ and assume $B \in \text{
pos}(\bar B,\cA)$ be given.  

We choose by induction on $n < \omega$ a sequence $\langle \bar
B_\eta:\eta \in {}^n 2 \rangle$ such that
\mn
\begin{enumerate}
\item[$\bullet$]  $\bar B_\eta \in \cB'$ moreover $(\exists^\infty
n)(B_{\eta,n} \backslash B_{\eta,n+1} \in \text{ id}^+_{\cA})$ 
for $\eta \in {}^n 2$
\sn
\item[$\bullet$]  $\bar B_\eta = \bar B$ if $\eta = \langle \rangle$
so $n=0$
\sn
\item[$\bullet$]  $B \in \text{ pos}(\bar B_\eta,\cA)$ if $\eta \in
{}^n 2$
\sn
\item[$\bullet$]  if $\nu \char 94 \langle 0 \rangle,\nu \char 94
\langle 1 \rangle \in {}^n 2$ then for some $k < \omega$ the set 
$B_{\nu \char 94 \langle 0 \rangle,k} \cap B_{\nu \char 94 \langle 1
\rangle,k}$ is finite.
\end{enumerate}
\mn
For $n=0$ this is trivial and for $n=m+1$ we use $\boxplus(c)$, i.e. the
construction of $\cB'$.  For every $n < \omega,\varrho \in {}^n 2$ let
$B_\varrho = \cap\{B_{\eta \rest k,m}:k \le n,m \le n\}$.  So
$B_\varrho \in \text{ id}^+_{\cA}$ and $m < \ell g(\varrho) \Rightarrow
B_\varrho \subseteq B_{\varrho \rest m}$ and if $\varrho_1 \ne \varrho_2
\in {}^n 2$ then for some $k < \omega$, for every $\rho_1 \in 
{}^{n+k} 2, \rho_2 \in {}^{n+k} 2$ satisfying 
$\varrho_1 \trianglelefteq \rho_1,\varrho_2 \trianglelefteq \rho_2$ we have
$B_{\rho_1} \cap B_{\rho_2}$ is finite. 
Obviously $[\varrho \in {}^\omega 2 \Rightarrow (\forall n <
\omega)(\exists k < \omega)(B_{\varrho,n} \backslash B_{\varrho,k} \in
\text{ id}^+_{\cA})]$ hence for each $\varrho
\in {}^\omega 2$ there is $C_\varrho \in \text{ id}^+_{\cA}$ such that
$C_\varrho \subseteq^* B_{\varrho \rest n}$ for $n < \omega$.

[Why?  We try by induction on $k < \omega$ to choose $A_{\varrho,k},
A'_{\varrho,k} \in \text{ob}(\omega)$ 
such that $A'_{\varrho,k} \in {\cA},A_{\varrho,k} 
\subseteq A'_{\varrho,k}$ and $m < k \Rightarrow A'_{\varrho,k} 
\neq A'_{\varrho,m}$ and $A_{\varrho,k} \subseteq^* B_{\varrho \rest k}$.
Now first, if we succeed then we can find $C \in \text{\rm ob}(\omega)$
such that for every $n < \omega$ we have $C \cap A_n$ is infinite and
$C \setminus \cup\{A'_{\varrho,m}:m < n\} \subseteq B_{\varrho \rest k_n}$.
If there is an infinite $C' \subseteq C$ almost disjoint
to every member of $\cA$, then $C_\varrho = C'$ is as required.
If there is no such $C'$ then we can find pairwise distinct 
$A''_n \in \cA \setminus \{A'_{\varrho,m}:m < \omega\}$
such that $C \cap A''_n$ is infinite for every $n < \omega$. 
Clearly $A''_n \cap C \subseteq^* B_{\varrho \rest m}$ for every $n,m < \omega$
and there is an infinite $C_\varrho \subseteq C$
such that $C _\varrho \subseteq^* B_{\varrho \rest m}$
and $C_\varrho \cap A''_n$ is infinite for every $n,m < \omega$, so
$C_\varrho$ is as required.

Second, if $k < \omega$ and  
we cannot choose $A_{\varrho,k}$ then we can choose 
$C_\varrho \in {\rm ob}(\omega)$ such that $n < \omega \Rightarrow
C_{\varrho} \subseteq^* B_{\varrho \rest n}$ and $C_{\varrho} 
\cap A_{\varrho,m} = \emptyset$ for $m < k$, and $C_\varrho$ is as required, so
we are done.] 

So $\cP' = \{B_{\eta \rest k,m}:k,m < \omega\}$ is as required.  

So proving $\boxplus$ is enough.

Why does this statement hold?

Let $\bar f = \langle f_\alpha:\alpha < {\gb}\rangle$
be a sequence of members of ${}^\omega \omega$ witnessing
${\gb}$ and \wilog \, $f_\alpha \in {}^\omega \omega$ is
increasing and $\alpha < \beta < {\gb} \Rightarrow f_\alpha
<_{J^{\text{bd}}_\omega} f_\beta$.

For $\alpha < {\gb}$ let $C_\alpha := \cup\{B_n \cap 
[0,f_\alpha(n)):n < \omega\}$ so clearly
\mn
\begin{enumerate}
\item[$(*)_1$]   $(a) \quad 
\alpha < \beta \Rightarrow C_\alpha \subseteq^* C_\beta$
\sn
\item[${{}}$]  $(b) \quad \alpha < \gb \wedge n < \omega \Rightarrow
C_\alpha \subseteq^* B_n$.
\end{enumerate}
\mn
We choose  $\alpha_\varepsilon = 
\alpha(\varepsilon) < {\gb}$ by induction on $\varepsilon <
{\frak b}$, increasing with $\varepsilon$ as follows: for $\varepsilon
= 0$ let $\alpha_\varepsilon = \text{ min}\{\alpha < {\gb}:
C_\alpha$ is infinite$\}$, for $\varepsilon = \zeta + 1$ let
$\alpha_\varepsilon = \text{ min}\{\alpha < {\gb}:\alpha >
\alpha_\zeta$ and $C_\alpha \backslash C_{\alpha(\zeta)}$ is infinite$\}$
and for $\varepsilon$ limit let $\alpha_\varepsilon =
\cup\{\alpha_\zeta:\zeta < \varepsilon\}$.  By the choice of $\bar f$
every $\alpha_\varepsilon$ is well defined, see the proof of 
$\oplus_\alpha$ below.

So $\langle \alpha_\varepsilon:\varepsilon < {\gb}\rangle$ is
increasing continuous with limit ${\gb}$.  
For each $\delta \in S^{\frak b}_{\aleph_0}$ let 
$\langle \varepsilon(\delta,n):n <
\omega\rangle$ be increasing with limit $\delta$ and, lastly, let
$\bar B_\delta = \langle C_{\alpha(\delta)} \backslash \bigcup\limits_{m \le n}
C_{\alpha(\varepsilon(\delta,m))}:n < \omega\rangle$ so $B_{\delta,n}
= C_{\alpha(\delta)} \backslash \bigcup\limits_{m \le n} 
C_{\alpha(\varepsilon(\delta,m))}$
hence $B_{\delta,n+1} \subseteq B_{\delta,n}$ and $B_{\delta,n}
\backslash B_{\delta,n+1}$ is infinite by the choice of
$\alpha_{\varepsilon(\delta,n)+1}$.  Clearly $\bar B_\delta \in {\cB}$
(also follows from the proof below). 

Why is $\langle \bar B_\delta:\delta \in S^{\gb}_{\aleph_0}\rangle$ as
required in $\boxplus$?  Clauses (a) + (b) are obvious and clause (d)
is easy (as if $\delta_1 < \delta_2$ then for some $n$ we have $\delta_1 <
\alpha(\varepsilon(\delta_2,n))$ hence $B_{\delta_1,n} \cap B_{\delta_2,n}
\subseteq B_{\alpha(\varepsilon(\delta_1,n))} \cap (B_{\alpha(\delta_2)} 
\backslash C_{\alpha(\varepsilon(\delta_2,n))})
\subseteq^* B_{\alpha(\delta_1)} \cap (B_{\alpha(\delta_2)} 
\backslash B_{\alpha(\delta_1)}) = \emptyset$.  
Lastly, to check clause (c) of $\boxplus$ let ${\cA}$
be an AD family and $B \subseteq \omega$ be such that

\mn
\begin{enumerate}
\item[$(*)_2$]
 $u = u_B := \{n < \omega:B \cap B_n \notin \text{ id}_{\cA}\}$ is
infinite, equivalently is $\omega$.
\end{enumerate}
\mn
It is enough to prove that for every $\alpha < {\gb}$
\mn
\begin{enumerate}
\item[$\oplus_\alpha$]    there is $\beta \in
(\alpha,{\frak b})$ such that $B \cap C_\beta \backslash C_\alpha
\in \text{ id}^+_{\cA}$.
\end{enumerate}
\mn
[Why is it enough?  As then for some club $E$ of ${\gb}$, for
every $\delta \in E \cap S^{\gb}_{\aleph_0}$ 
we have $(\forall \varepsilon < \delta)(\alpha_\varepsilon < \delta)$
and $(\forall \alpha < \delta)(\exists \beta)(\alpha < \beta < 
\delta \wedge C_\beta \backslash
C_\alpha \in \text{ id}^+_{\cA})$ hence $(\exists^\infty n)
((C_{\alpha(\varepsilon(\delta,n+1))} \backslash
C_{\alpha(\varepsilon(\delta,n))}) \in \text{ id}^+_{\cA})$ which means
$(\exists^\infty n)(B_{\delta,n} \backslash B_{\delta,n+1}) \in \text{
id}^+_{\cA})$ as required.]

So let us prove $\oplus_\alpha$.

If $\oplus_\alpha$ fails, for every $\beta \in (\alpha,{\gb})$ there are
$n=n(\beta)$ and $A_{\beta,0},\dotsc,A_{n(\beta)-1} \in {\cA}$ such that
$B \cap C_\beta \backslash C_\alpha \subseteq^* A_{\beta,0} \cup \ldots \cup
A_{\beta,n(\beta)-1}$.  Without loss of generality $n(\beta)$ is
minimal hence by $(*)_1$ the sequence 
$\langle n(\beta):\beta \in [\alpha,{\frak b})\rangle$ is
non-decreasing, but ${\gb} = \text{ cf}({\gb}) > \aleph_0$, 
hence, for some $\alpha_* \in [\alpha,{\gb})$, the
sequence $\langle n(\beta):\beta \in [\alpha_*,{\gb})\rangle$ is
constant and let $n(\alpha_*) = n_*$.

As ${\cA}$ is AD and $B \cap C_{\alpha_*} \backslash C_\alpha \subseteq^*
A_{\alpha_*,0} \cup \ldots \cup A_{\alpha_*,n_*-1}$ and $\beta \in
(\alpha_*,{\frak b}) \Rightarrow B \cap C_{\alpha_*} \backslash C_\alpha
\subseteq B \cap C_\beta \backslash C_\alpha \subseteq A_{\beta,0} \cup
\ldots \cup A_{\beta,n_*-1}$, using ``${\cA}$ is almost disjoint" and
the minimality of $n_{\alpha_*} = n_*$ 
it follows that $\{A_{\alpha_*,\ell}:\ell < n_*\} 
\subseteq \{A_{\beta,\ell}:\ell < n_*\}$ hence they are equal.

So
\mn
\begin{enumerate}
\item[$\odot$]   $\beta \in (\alpha,{\gb}) \Rightarrow B \cap 
C_\beta \backslash
C_\alpha \subseteq^* A_{\alpha_*,0}  \cup \ldots \cup
A_{\alpha_*,n_*-1}$.
\end{enumerate}
\mn
For each $n \in u = u_B$ as $B \cap B_n \backslash C_\alpha 
\in \text{ id}^+_{\cA}$, and
$A_{\alpha_*,0},\dotsc,A_{\alpha_*,n_*-1}$ are from 
id$_{\cA}$, clearly there is $k_n \in (B \cap B_n \backslash C_\alpha) 
\backslash A_{\alpha_*,0} \backslash \ldots \backslash
A_{\alpha_*,n_*-1} \backslash \{k_0,\dotsc,k_{n-1}\}$.
By the choice of $\bar f$ there is
$\beta \in (\alpha_*,{\frak b})$ such that $u_1 := \{n < \omega:k_n <
f_\beta(n)\}$ is infinite.  As $f_\beta$ is increasing, clearly $n \in u_1
\Rightarrow k_n < f_\beta(n) \Rightarrow k_n \in C_\beta \backslash
C_\alpha$.
So $\{k_n:n \in u_1\} \in [\omega]^{\aleph_0}$ is infinite and is 
a subset of $B \cap C_\beta
\backslash C_\alpha \backslash A_{\alpha_*,0},
\dotsc,A_{\alpha_*,n_*-1}$, so $\oplus_\alpha$ indeed holds, so we are done.
\end{PROOF}

\begin{PROOF}{\ref{g.3}}
\underline{Proof of \ref{g.3}}  We prove part (2), and part (1) follows
from it.  We immitate the proof of \ref{4d.3}.
\bigskip

\noindent
\underline{Stage A}:

Let $\kappa = {\gs}$.  Let ${\cP} \subseteq
[\kappa]^{\aleph_0}$ witness $\bold U(\kappa) = \kappa$, 
for transparency we assume $\omega \in \cP$ and 
$u \in {\cP} \Rightarrow \text{ otp}(u) = 
\omega$, this holds \wilog \, as ${\gb} \le \ga_* = \gs =  \kappa$.

\noindent
[Why?  It is enough to show that for every countable
$u \subseteq \kappa$ there is a family ${\cP}_u$ 
of cardinality $\le {\frak b}$ of subsets of $u$ 
each of order type $\omega$ such that every infinite
subset of $u$ has an infinite intersection with some member of $\cP$.  
Without loss of generality 
$u$ is a countable ordinal $\alpha$ and we prove this
by induction on $\alpha$.   For $\alpha$ successor ordinal or not 
divisible by $\omega^2$ this is trivial so let $\langle \alpha_n:
n < \omega \rangle$ be an increasing sequence of limit ordinals 
with limit $\alpha$ but $\alpha_0=0$. Let $\langle \beta_{n,k}:k 
< \omega \rangle$ list $[\alpha_n,\alpha_{n+1})$ with no repetitions and
let $\langle f_\epsilon \in {}^{\omega}\omega:\epsilon < {\frak b} 
\rangle$ exemplifies ${\frak b}$, each $f_\epsilon$ increasing and let
${\cP}_\alpha = \cup\{{\cP}_\beta:\beta < \alpha\} \cup 
\{\{\beta_{n,k}:n < \omega,k < f_\epsilon (n)\}:\varepsilon < \gb\}$. 
Clearly ${\cP}_\alpha$ has the right form and cardinality. 

Lastly, assume $v \subseteq u$ is infinite, if for some $\gamma <
\alpha,u \cap \gamma$ is infinite use the choice of $\cP_\gamma$. 
Otherwise let $f \in {}^{\omega}\omega$ be defined by 
$f(n) = \min \{k:(\exists m)[n \le m \wedge \beta_{m,k} \in v]\}$, and use
$\epsilon < {\frak b}$ large enough.]

Let $\langle u_\alpha:\alpha < \kappa\rangle$ list ${\cP}$ possibly
with repetitions, \wilog \, $n \le \omega \Rightarrow u_n = \omega$ 
and $\alpha > \omega \Rightarrow u_\alpha 
\subseteq \alpha$.  For $\alpha < \kappa$
let $\langle \gamma(\alpha,k):k < \omega\rangle$ list $u_\alpha$
in increasing order and $\gamma_{\alpha,k} = \gamma(\alpha,k)$.

Let $\langle {\cU}_\alpha:\alpha < \kappa\rangle$ be a partition of
$\kappa$, to sets each of cardinality $\kappa$ 
such that min$({\cU}_{1 + \alpha}) \ge \sup(u_\alpha) +1$ and $\omega
\subseteq {\cU}_0$. 
Let $\langle C^*_\alpha:\alpha \in {\cU}_0\rangle$ list a subset of 
${\cP}(\omega)$ witnessing ${\gs} = \kappa$ and as in Stage A of the
proof of \ref{4d.3}, the set $J_* \subseteq \text{ ob}(\omega)$ is
dense and $\omega \notin \text{ id}_{J_*}$.

If $\bar B$ is as in the assumption of \ref{g.19} and $\alpha \in
(0,\kappa)$ let ${\cP}_{\bar B}$ be as in the conclusion of \ref{g.19}
and for $\alpha < \kappa$ let 
$\langle C^*_{\bar B,\alpha,i}:i \in {\cU}_\alpha\rangle$ 
list ${\cP}_{\bar B}$.
\bigskip

\noindent
\underline{Stage B}: As in the proof of \ref{4d.3} but we use
$C^s_\rho(\rho \in {\cT}_s)$ which may really depend on $s$
and where $C^t_\rho,\bar B^t_{\nu,\beta}$ 
are defined in clauses $\boxplus_1(e),(g),(h),(i),(j)$ below 
(so the $\boxplus(e),(g),(h)$ from \ref{4d.3} are replaced) and depend just on 
${\cT}_t,\bar A_t$ and $\bar I_t$, too\footnote{also here we
require $\eta \in \text{ suc}({\cT}_s) \Rightarrow A_\eta \ne
\emptyset$} where 
\mn
\begin{enumerate}
\item[$\boxplus_1$]   (a)-(d) and (f) as in \ref{4d.3} of course and
\sn
\item[${{}}$]  $(e) \quad \bullet_1 \quad$ as before, i.e. $\bar A =
\langle A^t_\eta:\eta \in \text{ suc}(\cT_t)\rangle$,
\sn
\item[${{}}$]  $\,\,\,\,\,\,\, \quad \bullet_2 \quad \bar C = \bar C_t =
\langle C^t_\eta:\eta \in \Lambda_t \rangle$ where $\Lambda_t =
\{\eta:\eta \in {}^i 2$ and $i \in \cU_0$

\hskip35pt  or $\alpha > 0,i \in
\cU_\alpha$ and $\eta \rest (\sup(u_\alpha)) \in \cT_t$ \underline{or} (for
\ref{g.21}) $\eta \in \cT_t\}$
\sn
\item[${{}}$]  $\,\,\,\,\,\,\, \quad \bullet_3 \quad$ we stipulate
$A_\eta = \emptyset$ if $\eta \in \cT \backslash \text{ suc}(\cT)$
\sn
\item[${{}}$]  $(g) \quad$ as in \ref{4d.3} but replacing
$C^*_{\eta \rest i}$ by $C^t_{\eta \rest i}$
\sn
\item[${{}}$]  $(h) \quad$ if $i \in \cU_0$ and $\nu \in
\cT_t \cap {}^i 2$ then $C^s_\nu = C^*_i$
\sn
\item[${{}}$]  $(i) \quad$ if $\beta \in (0,\kappa)$ and $\nu
\in {}^{\sup(u_\beta)}2$ and both $\langle C^t_{\nu \rest i}:i \in
u_\beta\rangle$ and

\hskip25pt $\langle A^t_{\nu \rest i}:i \in u_\beta\rangle$ are well
defined \then \, we let $\bar B^t_{\nu,\beta} =
\langle B^t_{\nu,\beta,n}:n < \omega\rangle$

\hskip25pt  be defined by $B^t_{\nu,\beta,n} = \cap\{(C^t_{\nu \rest 
\gamma(\beta,k)})^{[\nu(\gamma(\beta,k)]} 
\backslash A^t_{\nu \rest \gamma(\beta,k)}:k < n\}$
\sn
\item[${{}}$]   $(j) \quad$ if $\beta \in (0,\kappa),i \in 
{\cU}_\beta$ hence $i \ge \sup(u_\beta)$ and 
$\rho \in {}^i 2$ and $\bar B^t_{\rho \rest 
\sup(u_\beta),\beta}$ 

\hskip25pt is well defined then, recalling stage A, $C^t_\rho = 
C^*_{\bar B^t_{\rho \rest \sup(u_\beta),\beta},\beta,i}$.
\end{enumerate}
\mn
Note that $\cT_t,\bar A_t$ determine $t$, i.e. $\bar
I_t,\Lambda_t,\bar C_t$ and $\langle \bar B^t_{\nu,\beta}:\nu,\beta$ as
above$\rangle$. 
\bigskip

\noindent
\underline{Stage C}:

As in \ref{4d.3} we just add:   
\mn
\begin{enumerate}
\item[$\boxplus_4$]   $(h) \quad$ if $s \le_{\text{AP}} t$ and $\bar
B^s_{\nu,\beta}$ is well defined then $\bar B^t_{\nu,\beta}$ is well
defined and

\hskip25pt  equal to it
\sn
\item[${{}}$]   $(i) \quad$ if $s \le_{\text{AP}} t$ and $C^s_\nu$
is well defined then $C^t_\nu$ is well defined and 

\hskip25pt  equal to it, so $\Lambda_s \subseteq \Lambda_t$
\sn
\item[${{}}$]   $(j) \quad C^s_\nu$ is well defined when 
$\nu \in c \ell({\cT}_s)  \cap {}^{\kappa >}2$.
\end{enumerate}
\mn
In the proof of $\boxplus_4(e)$ use the choice of $\langle C^s_\nu:\nu
\in {}^i 2,i \in {\cU}_0\rangle$, i.e. of 
$\langle C^*_\alpha:\alpha \in \cU_0\rangle$ in Stage A.

\bigskip

\noindent
\underline{Stages D,E}:  As in \ref{4d.3}.
\bigskip

\noindent
\underline{Stage F}:  The only difference is in the proof of $\odot_2$ in
subcase(2B).  Recall
\bigskip

\noindent
\underline{Case 2}:  SP$_B = \emptyset$ but not Case 1, i.e. $S_B
\subseteq \cT_s$, recall $B \subseteq \omega,B \notin J_s$.
\bigskip

\noindent
\underline{Subcase 2B}:  $\nu^*_B \notin S_B$ where $\nu^*_B =
\cup\{\eta:\eta \in S_B\}$
\mn
\begin{enumerate}
\item[$\odot_2$]   there is a set $B_1$ such that
\begin{enumerate}
\item[$(a)$]  $B_1 \subseteq B$ is infinite
\sn
\item[$(b)$]  $B_1$ is almost disjoint to any $A \in {\cA}^*$
\sn
\item[$(c)$]  if $\Lambda$ is finite then $\nu \in \Lambda 
\Rightarrow |B_1 \cap A_\nu| < {\aleph_0}$
\sn
\item[$(d)$]  if $\Lambda$ is infinite then for infinitely
many $\nu \in \Lambda$ we have $|B_1 \cap A_\nu| = {\aleph_0}$.
\end{enumerate}
\end{enumerate}
\mn
Why $\odot_2$ holds?  If $|{\cA}^*| < \kappa$ then ${\cA}^*$
has cardinality $< \kappa = {\gs}$ hence by the theorem's assumption
$|{\cA}^*| < {\frak s} = {\ga}_*$, so $\odot_2$ follows
as in the proof of \ref{4d.3}.  So we can assume $|{\cA}^*| 
= \kappa$ but $|{\cA}^*| \le \aleph_0 + |\ell g(\nu^*_B)|$ 
hence necessarily $\ell g(\nu^*_B) = \kappa$ follows and let 

\begin{equation*}
\begin{array}{clcr}
\cW := \{\alpha < \kappa:&\text{ for some } \ell \le n(\alpha) \text{ we
have } A^*_{\alpha,\ell} \cap B \in {\cA}^* \\
  &\text{ (equivalently } A^*_{\alpha,\ell} \cap B \text{ is infinite) but} \\
  &A^*_{\alpha,\ell} \notin \{A^*_{\alpha_1,\ell_1}:\alpha_1 < \alpha
\text{ and } \ell_1 \le n(\alpha_1)\}\}.
\end{array}
\end{equation*}

For $\alpha \in \cW$ choose $\ell(\alpha) \le n(\alpha)$ such that
$B \cap A^*_{\alpha,\ell(\alpha)}$ is infinite and $A^*_{\alpha,\ell(\alpha)}
\notin \{A^*_{\alpha_1,\ell_1}:\alpha_1 < \alpha$ and $\ell_1 \le
n(\alpha_1)\}$, in fact by $\odot_1(b)$ the last condition follows.  
As $n(\alpha) < \omega$ for $\alpha < \kappa$,
clearly $|\cW| = \kappa$ because $|{\cA}^*| = \kappa$, 
hence by the choice of ${\cP}$ there is
$u_* \in {\cP}$ such that $|\cW \cap u_*|$ is infinite; let 
$\alpha(*) \in [\omega,\kappa)$ be such that $u_{\alpha(*)} =
u_*$ and let $\nu = \nu^*_B \rest \sup(u_*)$; recall that otp$(u_*) =
\omega$; note that
\mn
\begin{enumerate}
\item[$\odot_{2.1}$]  $k <  \omega \Rightarrow 
B^s_{\nu,\alpha(*),k+1} \subseteq
B^s_{\nu,\alpha(*),k} \subseteq \omega$.
\end{enumerate}
\mn
[Why?  By their choice in $\boxplus_1(i)$.]

Recall also that 
$\langle \gamma_{\alpha(*),k}:k < \omega\rangle$ list $u_*$ 
in increasing order and so
\mn
\begin{enumerate}
\item[$\odot_{2.2}$]   $v := \{k < \omega:\gamma_{\alpha(*),k} 
\in \cW\}$ is infinite;
\sn
\item[$\odot_{2.3}$]   $[B^s_{\nu,\alpha(*),k}]^{\aleph_0} 
\supseteq I^s_{\nu \rest \gamma(\alpha(*),k)}$ for $k < \omega$.
\end{enumerate}
\mn
[Why?  As for $k(1) < k,(C^s_{\nu \rest
\gamma(\alpha(*),k(1))})^{[\nu(\gamma(\alpha)(*),k(1))]}$ and
$\omega \backslash A^*_{\nu \rest \gamma(\alpha(*),k(1))}$ belongs to
$\{X \subseteq \omega:[X]^{\aleph_0} \supseteq I^s_{\nu \rest
\gamma(\alpha(*),k)}\}$ hence by the definition of
$B^s_{\nu,\alpha(*),k}$ in $\boxplus_1(i)$ it satisfies $\odot_{2.3}$.]
\mn
\begin{enumerate}
\item[$\odot_{2.4}$]  $k \in v \Rightarrow B \cap 
B^s_{\nu,\alpha(*),k} \backslash
B^s_{\nu,\alpha(*),k+1}$ is infinite. 
\end{enumerate}
\mn
[Why?  For $k \in v$ let $\beta = \gamma(\alpha(*),k),n=n(\beta)$ and $\ell =
\ell(\beta)$.  On the one hand $[B \cap A^*_{\beta,\ell}]^{\aleph_0}
\subseteq [A^*_{\beta,\ell}]^{\aleph_0} \subseteq I_{\nu^*_B \rest
\beta}$.  On the other hand $A^*_{\beta,\ell}$ is disjoint to
$(C^s_{\nu^*_B \rest \beta})^{[\nu^*_B(\beta)]} \backslash A^s_{\nu^*_B \rest
\beta}$ if $A^*_{\beta,\ell} = A^*_{\beta,n}$ trivially and 
is almost disjoint to $(C^s_{\nu^*_B \rest
\beta})^{[\nu(\beta)]}$ otherwise (i.e. as $[(C^s_{\nu^*_B \rest
\beta})^{[1-\nu(\beta)]}]^{\aleph_0} \supseteq I_{(\nu^*_\beta \rest
\beta) \char 94 \langle 1 - \nu(\beta)\rangle} \supseteq
\{A^*_{\alpha,\ell}\}$).  Hence 
$(C^s_{\nu^*_B \rest \beta})^{[\nu^*_B(\beta)]} 
\backslash A^*_{\beta,\ell}$ is almost disjoint to
$B \cap A^*_{\beta,\ell}$,  an infinite set from 
$I_{\nu^*_B \rest \beta}$ hence by $\odot_{2.3}$ from
$[B^s_{\nu,\alpha(*),k}]^{\aleph_0}$.  So $B \cap B^s_{\nu,\alpha(*),k}
\backslash B^s_{\nu,\alpha(*),k+1} = B \cap B^s_{\nu,\alpha(*),k} \backslash
(B^s_{\nu,\alpha(*),k} \cap ((C^s_{\nu^*_B \rest
\beta})^{[\nu^*_B(\beta)]} \backslash A^*_{\beta,k})$ almost contains
this infinite set hence is infinite as promised.]

So by the choice of ${\cP}_{\bar B^s_\nu,\alpha(*)}$, i.e. 
\ref{g.19} and clauses (i),(j) of
$\boxplus_1$ for some $\beta \in {\cU}_{\alpha(*)}$ so $\beta \ge
\alpha(*) \ge \ell g(\nu)$ we have $B \backslash
(C^s_{\nu^*_B \rest \beta})^{[\ell]} \notin J_s$ for $\ell=0,1$
hence $B_1 := B \backslash (C^s_{\nu^*_B})^{[\nu^*_B(\beta)]}
\notin J_s$ recalling that for $\beta \in {\cU}_\alpha,\alpha
\ne 0$ and $\rho \in {}^\beta 2$ the set $C^s_\rho$ depends just on
$\ell g(\rho)$ and $\rho \rest \sup(u_\alpha)$ (and our $s$).

Now consider $B_1$ instead of $B$, clearly $S_{B_1}$ is a subset of
$S_B$ and $\nu^*_B \rest (\beta +1)$ is not in it, but $B_1 \in \text{
ob}(B)$ hence $S_{B_1} \subseteq S_B$ hence $S_{B_1}$ is $\subseteq
\{\nu^*_B \rest \gamma:\gamma \le \beta\}$ and $B_1$ fall under subcase
(2A) as $\beta < \kappa = \ell g(\nu^*_B)$.
\end{PROOF}

\begin{theorem}
\label{g.21} 
There is a saturated MAD
family ${\cA} \subseteq J_*$ when ${\ga}_* < \kappa = {\gs},J_*
\subseteq \text{\rm ob}(\omega)$ is dense and 
{\rm Pr}$(\kappa,{\ga})$, see \ref{g.5}(3). 
\end{theorem}

\begin{PROOF}{\ref{g.21}}
\underline{Proof of \ref{g.21}}  
We immitate the proofs of \ref{4d.3}, \ref{g.3}.  
Note that ${\gb} \le {\ga}_* < {\gs}$.
\bigskip

\noindent
\underline{Stage A}:

Similarly to stage A of the proof of \ref{g.3}; let $(E,\bar{\cP}^*)$
be as in Definition \ref{g.5}(3) and Observation \ref{g.10}, 
as ${\gb} < \kappa$, \wilog \, $u \in {\cP}^*_\alpha \Rightarrow \text{ otp}(u)
= \omega$ for $\alpha \in [\omega,\kappa)$.  
As we can replace $E$ by any appropriate club $E'$ of $\kappa$ contained in
acc$(E)$ see \ref{g.10} there, \wilog \, otp$(E) = \text{
cf}(\kappa)$, min$(E) \ge \omega$ and $\gamma \in E \Rightarrow \gamma + 1
+ \gb < \text{ min}(E \backslash (\gamma +1))$.  
Let $\langle \gamma^*_i:i < \text{ cf}(\kappa)\rangle$ 
list $E$ in increasing order.

Let $\langle u_\gamma:\gamma < \kappa\rangle$ be such that 
$\langle u_\gamma:\gamma^*_i \le \gamma <
\gamma^*_{i+1}\rangle$ list ${\cP}_{\gamma^*_{i+1}}$ (which 
includes $\cP_{\gamma_i}$) 
and $u_j = \omega$ for $j < \gamma^*_0$.

Let $\langle {\cU}_\alpha:\alpha < \kappa\rangle$ be a partition of
$\{2i +1:i < \kappa\}$ such that 
min$({\cU}_{1 + \alpha}) \ge \alpha + \omega,|{\cU}_{1 + \alpha}| = 
{\gb},|{\cU}_0|=\kappa,1 \le \alpha < \gamma^*_i \Rightarrow
{\cU}_\alpha \subseteq \gamma^*_i$.

Let $\langle C^*_i:i \in {\cU}_0\rangle$ list a family of subsets
of $\omega$ witnessing ${\gs} = \kappa$ also $J_*$ is as in the proof
of \ref{4d.3}.

Let ${\cP}_{\bar B},\langle C^*_{\bar B,\alpha,i}:i \in
{\cU}_\alpha\rangle$ be as in \ref{g.3}, Stage A.
\bigskip

\noindent
\underline{Stage B}:

As in \ref{g.3}, i.e. the case ${\gs} = {\ga}$, but we
change $\boxplus_1(f)$
\mn
\begin{enumerate}
\item[$\boxplus_1$]  $(f) \quad \bullet \quad A_\eta \in I_\eta \cap J_*$ or
$A_\eta = \emptyset$ and 
\sn
\item[${{}}$]  $\qquad \,\,\,\,\bullet \quad \mathscr{S}_t := 
\{\eta \in \cT_t:A_\eta \ne \emptyset\}
\subseteq \{\eta \in {\cT}_t:\ell g(\eta) = \gamma^*_i +1$

\hskip25pt  for some $i < \kappa\}$.
\end{enumerate}
\bigskip

\noindent
\underline{Stage C}:  

As in the proof of \ref{g.3}.
\bigskip

\noindent
\underline{Stage D}:

Here there is a minor change: we replace $\boxplus_7$ in \ref{4d.3},
\ref{g.3} by $\boxplus_7,\boxplus_8,\boxplus_9$ below
\mn
\begin{enumerate}
\item[$\boxplus_7$]  if $\alpha < 2^{\aleph_0},s \in 
\text{ AP}_\alpha$ and $B \in J^+_s$ \then \, there are a limit
ordinal $\xi \in \kappa \backslash E$ and 
$t \in \text{ AP}_{\alpha +1}$ such that $s \le_{\text{AP}} t$ 
and $|S^t_B \cap {}^\xi 2| = 2^{\aleph_0}$; we may add $\mathscr{S}_t
= \mathscr{S}_s$. 
\end{enumerate}
\mn
This is proved in Stage F.

To clarify why this is O.K. recall $\boxplus_6(f)$ and note that
\mn
\begin{enumerate}
\item[$(*)$]  if $s \le_{\text{AP}} t,B \in \text{ ob}(\omega),\eta
\in S^s_B \backslash \cT_s$ and $\eta \notin \cT_t$ then
$\eta \in S^t_B$.
\end{enumerate}
\mn
Now we need
\mn
\begin{enumerate}
\item[$\boxplus_8$]   if $\xi \in \kappa \backslash E$ is a limit
ordinal, $\alpha < 2^{\aleph_0},t \in \text{ AP}_\alpha,B \in \text{
ob}(\omega)$ and $|S^t_B \cap {}^\xi 2| = 2^{\aleph_0}$ and $\zeta
= \text{ min}(E \backslash \xi)$ \underline{then} for every $t_1$
and $\alpha + \zeta \le \beta < 2^{\aleph_0}$
 such that $t \le_{\text{AP}}
t_1 \in \text{ AP}_\beta$ there is $t_2,t_1 \le_{\text{AP}} t_2 \in
\text{ AP}_{\beta +1}$ satisfying $(\exists \eta \in \text{
suc}({\cT}_{t_2}))[\eta \notin {\cT}_{t_1} \wedge A^{t_2}_\eta \in
\text{ ob}(\omega) \wedge A^{t_2}_\eta \subseteq B]$.
\end{enumerate}
\mn
The proof of $\boxplus_8$ is like the proof of Case 1 in Stage F in
the proof of \ref{4d.3} but we elaborate; we are given
$\beta,\xi,\zeta$ and $t_1$ such
that $t \le_{\text{AP}} t_1 \in \text{ AP}_\beta$; now we choose
$\rho \in S^t_B \cap {}^\xi 2 \backslash {\cT}_{t_1}$ exists as 
$|S^t_B  \cap {}^\xi 2| = 2^{\aleph_0} > |{\cT}_{t_1}|$ so recalling
$(*)_5(e)$ necessarily $\rho \in S^{t_1,1}_B$.
Choose $B_1$ such that $B_1 \subseteq B,B_1 \in I^t_\rho$.

Note that for every $\varepsilon \in [\xi,\zeta +1)$ either
$C^{t_1}_\varrho$ is well defined for every $\varrho \in
{}^\varepsilon 2$ such that $\rho \trianglelefteq \varrho$ and its
value is the same for all such $\varrho$ (when $\varepsilon$ is odd) 
\underline{or} $C^{t_1}_\varrho$ for $\rho
\trianglelefteq \varrho \in {}^\varepsilon 2$ is not well defined
(when $\varepsilon$ is even).  So
${\cB} = \{C^{t_1}_\varrho:\rho \triangleleft \varrho \in {}^{\zeta +1
\ge}2$ and $C^{t_1}_\varrho$ is well defined$\}$ is a family of 
$\le |\zeta| < \kappa = {\frak s}$ subsets of
$B_1$ hence there is an infinite $B_2 \subseteq B_1$ such that $\rho
\trianglelefteq \varrho \in {}^{\zeta >}2 \wedge (C^t_\varrho$ well
defined) $\Rightarrow B_2 \subseteq^* C^t_\varrho \vee B_2 \subseteq^*
\omega \backslash C^t_\varrho$ and \wilog \,
$B_2 \in J_*$.

We choose $\eta$ such that $\rho \triangleleft \eta \in {}^{\zeta
+1}2$ and 

$[\ell g(\rho) \le \gamma < \zeta +1 \wedge (C^t_{\eta
\rest \gamma}$ is well defined) $\wedge B_2 \subseteq^* (C^t_{\eta
\rest \gamma})^{[\ell]} \wedge \rho \in \{0,1\} \Rightarrow
\eta(\gamma) = \ell]$.  
Let us define $t_2 \in \text{ AP}_{\beta + \zeta +2} := 
\text{ AP}_{\beta +1}$ (as $\alpha + \zeta +1 \le \beta$ and
$|\alpha_1| = |\alpha_2| \Rightarrow \text{ AP}_{\alpha_1} = \text{
AP}_{\alpha_2}$) as follows:
\mn
\begin{enumerate}
\item[$(a)$]   ${\cT}_{t_2} := {\cT}_{t_1}
\cup\{\varrho:\varrho \trianglelefteq \eta\}$
\sn
\item[$(b)$]   $A^{t_2}_\varrho$ is $A^{t_1}_\varrho$ if well
defined, is $B_2$ if $\varrho = \eta$ and is $\emptyset$
if $\eta \in \text{ suc}({\cT}_{t_2})$ but $A^{t_2}_\varrho$ is not
already defined
\sn
\item[$(c)$]   $C^{t_2}_\varrho$ is $C^{t_1}_\varrho$ if $\varrho
\in {\cT}_{t_1}$ and we choose $C^{t_2}_{\eta \rest \varepsilon}$
by induction on $\varepsilon \in [\xi,\zeta +2]$ as follows: if it is
determined by $\boxplus_1$ we have no choice otherwise let it be
$\omega^{[\eta(\varepsilon)]}$.
\end{enumerate}
\mn
The other objects of $t_2$ are determined by those we have chosen.  So
$\boxplus_8$ holds indeed.
\mn
\begin{enumerate}
\item[$\boxplus_9$]   if $s \in \text{ AP}_\alpha$ and $\rho \in c
\ell({\cT}_s)$ then for some $t,s \le_{\text{AP}} t \in \text{
AP}_{\alpha + 3}$ and ${\cT}_s \subseteq {\cT}_t \subseteq 
{\cT}_s \cup \{\rho,\rho \char 94
\langle 0 \rangle,\rho \char 94 \langle 1\rangle\}$ and  $I^s_\rho \ne
\emptyset \Rightarrow \rho \in {\cT}_t$ and $I^s_\rho \ne \emptyset
\wedge \ell < 2 \wedge I^t_{\rho \char 94 <\ell>} \ne \emptyset
\Rightarrow \rho \char 94 \langle \ell \rangle \in {\cT}_t$.
\end{enumerate}
\mn
[Why?  Easier than $\boxplus_8$.]
\bigskip

\noindent
\underline{Stage E}:

Similar to \ref{4d.3} with the changes necessitated by the change in
Stage D.
\bigskip

\noindent
\underline{Stage F}:

We prove $\boxplus_7$, the proof splits to cases.
\bigskip

\noindent
\underline{Case 1}:  Some $\nu \in S_B$ is such that $\nu \in c 
\ell({\cT}_s) \backslash {\cT}_s$.

Let $B_1 \in \text{ ob}(B) \cap I^s_\nu$, there is such $B_1$ as $\nu
\in S_B$ but $\nu \notin S^{s,2}_B$ as $\nu \notin \cT_s$.

Let $C_{\nu,n} \in \text{ ob}(\omega)$ for $n < \omega$ be such that
$\cap\{C^{[\varrho(n)]}_{\nu,n}:n < \ell g(\varrho)\} \cap B_1$
is infinite for every $\varrho \in {}^{\omega >}2$.

We choose ${\cT}_t = T_s \cup \{\nu \char 94 \rho:\rho \in
{}^{\omega >}2\}$.  For $\rho \in {}^{\omega >}2$, we choose $C^t_{\nu
\char 94 \rho}$ by induction on $\ell g(\rho)$: if $\ell g(\nu \char
94 \rho) = \ell g(\nu) + n$ is even and $n \in \{2m,2m+1\}$ then
$C^t_{\nu \char 94 \rho} = C_{\nu,m}$, otherwise we act as in the
proof of $\boxplus_8$.  Lastly, let $A^t_{\nu \char 94 \rho} =
\emptyset$ for $\rho \in {}^{\omega >}2$.

Easily
\mn
\begin{enumerate}
\item[$(*)$]    if $s \le_{\text{AP}} s_1,|{\cT}_{s_1}| < 2^{\aleph_0}$
then $|S^{s_1}_B \cap {}^{\ell g(\nu)+\omega}2| = 2^{\aleph_0}$.
\end{enumerate}
\mn
So we are done with Case 1.
\bigskip

\noindent
\underline{Case 2}:  $SP^s_B = \emptyset$ but not Case 1 and let $\nu^*_B =
\cup\{\eta:\eta \in S_B\}$.
\bigskip

\noindent
\underline{Subcase 2A}:  $\nu^*_B \in S_B$

As in the proof \ref{g.3} but end as in Case 1.
\bigskip

\noindent
\underline{Subcase 2B}:  $\nu^*_B \notin S_B$

Except $\odot_2$ which we elaborate
this is as in the proofs of \ref{4d.3}, \ref{g.3} 
but in the end replace ``Subcase 1B" by ``Case 1".  Recall from 
Stage 2B as in the proof of \ref{4d.3}, $\delta := \ell g(\nu^*_B)$
is a limit ordinal and $\nu^*_{B,\alpha} = (\nu^*_B \rest \alpha) \char 94
\langle 1 - \nu^*_B(\alpha)\rangle$ for $\alpha < \delta$  
and $\langle \langle A^*_{\alpha,n}:n \le n(\alpha)\rangle:
\alpha < \delta\rangle$ and $\cA^*
:= \{B \cap A:A = A^*_{\alpha,k}$ for some $\alpha < \delta,k \le
n(\alpha)$ and $B \cap A$ is infinite$\}$, $\Lambda = \{\nu \in
\mathscr{S}_s:\nu^*_B \trianglelefteq \nu$ and $A_\nu \cap B$ is
infinite$\}$ are as in Stage (2B) of the
proof of \ref{4d.3}.  We have to prove:
\mn
\begin{enumerate}
\item[$\odot_2$]   there is a set $B_1$ such that:
\begin{enumerate}
\item[$(a)$]   $B_1 \subseteq B$ is infinite
\sn
\item[$(b)$]  $B_1$ is almost disjoint to any $A \in {\cA}^*$
\sn
\item[$(c)$]  if $\Lambda$ is finite then
$\nu \in \Lambda \Rightarrow |B_1 \cap A_\nu| < {\aleph_0}$
\sn
\item[$(d)$]  if $\Lambda$ is infinite then for infinitely
many $\nu \in \Lambda$ we have $|B_1 \cap A_\nu| = {\aleph_0}$.
\end{enumerate}
\end{enumerate}
\mn
Why $\odot_2$ holds?  If $|{\cA}^*| < \ga_*$ then as in 
the proof of \ref{4d.3} the statement of $\odot_2$ follows.
So we can assume $|{\cA}^*| \ge \ga_*$ and let 

\begin{equation*}
\begin{array}{clcr}
\cW := \{\alpha < \kappa:&\text{ for some } \ell \le n(\alpha) \text{ we
have } A^*_{\alpha,\ell} \cap B \in {\cA}^* \\
  &\text{ (equivalently } A^*_{\alpha,\ell} \cap B \text{ is infinite)}\}.
\end{array}
\end{equation*}
\mn
and let

\[
\cW' = \{\alpha \in W:|\alpha \cap \cW| < \ga_*\}.
\]
\mn
\underline{Subcase 2B$(\alpha)$}:  sup$(\cW') \in \text{ acc}(E')$

For $\alpha \in \cW$ choose $\ell(\alpha) \le n(\alpha)$ such that
$B \cap A^*_{\alpha,\ell(\alpha)}$ is infinite hence $A^*_{\alpha,\ell}
\notin \{A^*_{\alpha_1,\ell_1}:\alpha_1 < \alpha,\ell_1 \le
n(\alpha_1)\}$.  As $n(\alpha) < \omega$ for $\alpha < \kappa$,
clearly $|\cW| = |\cA^*| \ge \ga$ hence otp$(\cW') = \ga$.

So by Definition \ref{g.5}(3), i.e. the
choice of $(E,\bar\cP)$, there is a pair $(u_*,\gamma^*_j)$ as in
clause (e) there.  So $u_* \in \cP^*_{\gamma^*_j}$ hence 
$u_* = u_{\alpha(*)}$ for some
$\alpha(*) \in [\gamma^*_j,\gamma^*_{j+1}),\gamma^*_{j+1} < \sup(\cW')$ 
and $\cW \cap \sup(E \cap \gamma^*_j) = \cW' \cap 
\sup(E \cap \gamma^*_j)$ has cardinality $< \ga_*$
and let $\nu = \nu^*_B \rest \sup(u_*)$; recall otp$(u_*) = \omega$.

Recall also that $\langle \gamma_{\alpha(*),n}:n < \omega\rangle$ list $u_*$
in increasing order
and so $v := \{n < \omega:\gamma_{\alpha(*),n} \in \cW\}$ is infinite and
clearly $n \in v \Rightarrow B^s_{\nu,\alpha(*),n} \backslash
B^s_{\nu,\alpha(*),n+1}$ is infinite as in the proof of \ref{g.3}.  So 
by the choice of ${\cP}_{\bar B^s_\nu,\alpha(*)}$, 
i.e. \ref{g.19} and clauses (i),(j) of
$\boxplus_1$ for some $\beta \in {\cU}_{\alpha(*)}$ 
so $\beta \ge \ell g(\nu)$ we have $B \backslash
(C^s_{\nu^*_B \rest \beta})^{[\ell]} \notin J_s$ for $\ell=0,1$
hence $B_1 := B \backslash (C^s_{\nu^*_B})^{[\nu^*_B(\beta)]}
\notin J_s$ recalling that for $\beta \in {\cU}_\alpha,\alpha
\ne 0$ and $\rho \in {}^\beta 2$ the set $C^s_\rho$ depends just on
$\ell g(\rho)$ and $\rho \rest \sup(u_\alpha)$ (and our $s$).
\bigskip

\noindent
We finish as in the proof of \ref{g.3}.
\medskip

\noindent
\underline{Subcase 2B$(\beta)$}:   sup$(\cW') \in (\gamma^*_i,\gamma^*_j]$
where $j=i+1$ so $\gamma^*_i,\gamma^*_j \in E$ and let 
$\gamma^* = \sup(\cW')$.

Apply Definition \ref{g.5}(3), clause (f) to $\gamma^*_i,\gamma^*_j,\cW'
\backslash \gamma^*_j$ we
get $u = u_* \in \cP_{\gamma^*_j}$ so $u = u_{\alpha(*)} \subseteq
[\gamma^*_i,\gamma^*_j)$ for some $\alpha(*) \in (\gamma^*_i,\gamma^*_j)$.

Let $\beta = \sup(u_*),\nu = \nu^*_\beta \rest \beta$ so by
$\odot_{2.1} + \odot_{2.4}$ in stage 2B in the proof of \ref{g.3},
and by the choice of $\cP_{\bar B_{\nu,\alpha(*)}}$ there is 
a $\cQ \subseteq \cP_{\bar B^s_{\nu,\alpha(*)}}$ 
of cardinality $\aleph_0$ and $\Lambda \subseteq
{}^{\cQ} 2$ of cardinality $2^{\aleph_0}$ such that for every $\rho \in
\Lambda$ there is $B_\rho \in \text{ ob}(B) \cap J^+_s$ such
that $A \in \cQ \Rightarrow B_\rho \subseteq^* A^{[\rho(A)]}$ and $n < \omega
\Rightarrow B_{\eta,n} :
= B_n \backslash \cup\{C^{[\nu^*_B(\gamma_{\alpha(*),k})]}_{\nu^*_B
\rest \gamma_{\alpha(*),k}}:k < n\} \in I^s_{\nu^*_B \rest
\gamma_{\alpha(*),n}}$.

Clearly for some $v \subseteq \cU_{\alpha(*)} \subseteq
(\gamma^*_i,\gamma^*_j)$ of cardinality $\aleph_0,\nu
\triangleleft \rho \in \cT_s \wedge \ell g(\rho) \ge \sup(v)
\Rightarrow \{C^s_{\rho \rest \varepsilon}:\varepsilon \in v\} = \cQ$.

For $\eta \in \Lambda$ analyzing $S_{B_\eta}$ and 
recalling $\gamma^*_{j+1} < \gs$
clearly $S_{B_\eta} \cap \{\nu \in
\cT_s:\ell g(\nu) \le \sup(u_*)\}$ is $\{\nu^*_B \rest \gamma:\gamma <
\sup(u_*)\}$ and sup$(u_*) > \gamma^* + 1$, so there is no $\rho$ such that
$A^s_\rho$ is non-empty, $\rho \in \cT_s$ and $\sup(u_*) \le \ell g(\rho) 
< \gamma^*_j$, so $S_{B_\eta} \cap 
\{\nu \in \cT_s:\ell g(s) <\gamma^*_j\}$ does not
depend on $\langle A^s_\eta:\eta \in \text{ suc}(\cT_s),\ell g(\eta)
\ge \gamma^*_j\rangle$ so we can finish easily as in case 1.
\bigskip

\noindent
\underline{Case 3}:   As in the proof of \ref{g.3}.
\end{PROOF}
\newpage

\section {Further Discussion} 

The cardinal invariant ${\gs}$ plays here a major role, so the
claims depend on how ${\gs}$ and ${\frak a}_*$ are compared; when
${\gs} = {\ga}_*$ it is not clear whether the further assumption of
\ref{g.3}(2) may fail.  If ${\gs} > {\ga}_* > \aleph_1$, it is
not clear if the assumption of \ref{g.21} may fail.  Recall
\ref{4d.3}, dealing with $\gs < \ga_*$, the first case is proved 
ZFC, but the others need pcf assumptions.

All this does not
exclude the case ${\gs} = \aleph_{\omega +1},{\ga}_* = \aleph_1$
hence ${\gb} = \aleph_1$, as in \cite{Sh:668}.  Fulfilling the
promise from \S0 and the abstract.

\begin{claim}
\label{4q.0}  
1) If there is no inner model with a measurable cardinal (and even the
non-existence of much stronger statements) \then
  \, there is a saturated MAD family, $\cA$.

\noindent
2) Also if $\gs < \aleph_\omega$ there is one.

\noindent
3) Moreover, if $J_* \subseteq \text{\rm ob}(\omega)$ is dense 
\then \, we can demand $\cA \subseteq J_*$.
\end{claim}

\begin{PROOF}{\ref{4q.0}}
As Theorems \ref{4d.3}, \ref{g.3}, \ref{g.21} cover
them (using well known results).
\end{PROOF}

\noindent
We now remark on some further possibilities.
\begin{definition}
\label{4q.1}
1) We say ${\mathscr S} \subseteq
\text{ ob}(\omega)$ is ${\gs}$-free \underline{when}:
\mn
\begin{enumerate}
\item[$(a)$]   for every $A \in \text{ ob}(\omega)$ there is $B \in
\text{ ob}(A)$ such that $B$ induces an ultrafilter on ${\mathscr S}$;
i.e. $C \in {\mathscr S} \Rightarrow A \subseteq^* C \vee A \subseteq^*
(\omega \backslash C)$.
\end{enumerate}
\mn
1A) We say ${\mathscr S} \subseteq \text{ ob}(\omega)$ is ${\gs}$-free
in $I$ \underline{when} $I \in \text{ OB}$ and for every $A \in I$ there is
$B \in \text{ ob}(A)$ which induces an ultrafilter on $S$.

\noindent
2) We say ${\mathscr S} \subseteq \text{ ob}(\omega)$ is 
${\gs}$-richly free \when \, clause (a) and
\mn
\begin{enumerate}
\item[$(b)$]    if $A \in \text{ ob}(\omega)$ and the set
$\{D \cap {\mathscr S}:D$ an ultrafilter on $\omega$ containing
ob$(A)\}$ is infinite, \underline{then} it has cardinality continuum.
\end{enumerate}
\mn
3) We say ${\mathscr S} \subseteq \text{ ob}(\omega)$ is 
${\frak s}$-anti-free if no $B \in \text{ ob}(\omega)$ induces an
ultrafilter on ${\mathscr S}$.

\noindent
4) Let ${\gS}$ be $\{\kappa$: there is a $\subseteq$-increasing sequence
$\langle {\mathscr S}_i:i < \kappa\rangle$ of ${\gs}$-richly-free families
such that $\cup\{{\mathscr S}_i:i < \kappa\}$ is not ${\gs}$-free.

\noindent
5) Recall ${\gs} = \text{ min}\{|{\mathscr S}|:{\mathscr S} \subseteq
\text{ ob}(\omega)$ and no $B \in \text{ ob}(\omega)$ induces an
ultrafilter on ${\mathscr S}\}$.

\noindent
6) We say ch$_{\text{dim}}({\cB}) < \kappa$ when $(\exists \eta \in
{}^{\cB}2)(I = I_{{\cB},\eta}) \Rightarrow 
\text{ cf}(I,\subseteq) < \kappa$ recalling \ref{0z.6}(5).
\end{definition}

\begin{observation}
\label{4q.4}
1) If ${\mathscr S}$ is ${\gs}$-free and ${\mathscr S}' \subseteq 
{\mathscr S}$ then ${\mathscr S}'$ is ${\gs}$-free.

\noindent
2) If ${\mathscr S} \subseteq \text{ ob}(\omega)$ and $|{\mathscr S}| 
< {\gs}$ then ${\mathscr S}$ is ${\gs}$-free.

\noindent
3) If ${\mathscr S}_n \subseteq \text{ ob}(\omega)$ is ${\gs}$-free for
   $n < \omega$ then $\cup\{{\mathscr S}_n:n < \omega\}$ is ${\gs}$-free.

\noindent
4) ${\gs} \in {\gS}$.

\noindent
5) $\kappa \in {\gS}$ iff cf$(\kappa) \in {\gS}$.

\noindent
6) $\kappa \in {\gS} \Rightarrow \aleph_1 \le \kappa \le 2^{\aleph_0}$.

\noindent
7) In Definition of ${\gS}$ we can add ``$\cup\{{\mathscr S}_i:i <
\kappa\}$ is ${\gs}$-anti-free".

\noindent
8) cf$({\gs}) > \aleph_0$, in fact $\kappa \in {\gS}
\Rightarrow \text{ cf}(\kappa) > \aleph_0$.
\end{observation}

\begin{definition}
\label{4q.7}  
1) We say $A \in \text{ ob}(\omega)$ obeys 
$f \in {}^\omega \omega$ \when \,: for every $n_1 <  n_2$ 
from $A$ we have $f(n_1) < n_2$.

\noindent
2) Let $\bar f = \langle f_\alpha:\alpha < \delta\rangle$ be a
   sequence of members of ${}^\omega \omega$.
We say $\bar A = \langle A_\alpha:\alpha\in u\rangle$ obeys 
$\bar f$ \underline{when} $u \subseteq \delta$ and $A$ obeys $f_\alpha$ for
$\alpha \in A$.

\noindent
3) ${\ga}_{\bar f} = \text{ Min}\{|u|$: there are $B \in \text{
   ob}(\omega)$ and $\bar A = \langle A_\alpha:\alpha \in u\rangle$
   obeying $\bar f$ such that $\{A_\alpha \cap B:\alpha \in u\}$ is a
   MAD of $B\}$.
\end{definition}

\begin{remark} 
\label{4q.12}
1)  Also note that in \ref{4d.3}, \ref{g.3}, \ref{g.21} 
we can replace ${\frak s}$ by a smaller (or equal)
cardinal invariant ${\frak s}_{\text{\rm tree}}$,
the tree splitting number.

\noindent
2) Let ${\frak s}_{\text{\rm tree}}$ be the minimal $\kappa$ such that
there is a sequence 
$\bar C = \langle C_\eta:\eta \in {}^{\kappa >} 2 \rangle$ 
such that $C_\eta \in {\rm ob}(\omega)$ for $\eta \in {}^{\kappa >} 2$  
and there is no $\eta \in {}^\kappa 2$ and $A \in {\rm ob}(\omega)$
such that $\epsilon < k \Rightarrow A \subseteq^* C_\eta^{[\eta(\epsilon]}$.
Note that the minimal $\kappa$ for which
there is such sequence $\langle C_\eta:\eta \in {}^{\kappa >} 2 
\rangle$ has uncountable cofinality.

\noindent
3) Also in \ref{4d.3} we may weaken ${\frak s} < {\frak a}_*$ to
${\frak s} < {\frak a} \wedge {\frak s} \le {\frak a}_*$. 
\end{remark}


\end{document}